\documentclass[dvipdffmx]{amsart}
\usepackage{amsmath,amssymb}
\usepackage{mathrsfs}
\usepackage{amscd}
\usepackage{booktabs}
\usepackage{graphicx,xcolor}
\newtheorem{df}{Definition}
\newtheorem{prp}{Proposition}
\newtheorem{thm}{Theorem}

\newtheorem{rem}{Remark}

\newcommand{\relmiddle}[1]{\mathrel{}\middle#1\mathrel{}}
\numberwithin{equation}{section}
\numberwithin{figure}{section}
\numberwithin{df}{section}
\numberwithin{prp}{section}
\numberwithin{thm}{section}
\numberwithin{lem}{section}
\numberwithin{rem}{section}
\numberwithin{con}{section}
\numberwithin{ex}{section}
\begin{document}

\title[Perfectness of Kirillov-Reshetikhin Crystals]{Perfectness of Kirillov-Reshetikhin Crystals $B^{r,s}$ for types $E_{6}^{(1)}$ and $E_{7}^{(1)}$ with a minuscule node $r$}

\author{Toya Hiroshima}
\address{OCAMI, Osaka City University,
3-3-138 Sugimoto, Sumiyoshi-ku, Osaka 558-8585, Japan}
\email{toya-hiroshima@outlook.jp}


\date{}

\begin{abstract}
We prove the perfectness of Kirillov-Reshetikhin crystals $B^{r,s}$ for types $E_{6}^{(1)}$ and $E_{7}^{(1)}$ with $r$ being the minuscule node and $s\geq 1$ using the polytope model of KR crystals introduced by Jang~\cite{Jan21}.
\end{abstract}


\maketitle


\section{Introduction}

Kirillov-Reshetikhin (KR) modules are an important class of finite-dimensional representations for affine quantum group $U_{q}^{\prime}(\mathfrak{g})$ without the degree operator~\cite{CP94}.
We denote a KR module by $W^{r,s}$, where $r$ is a node of the Dynkin diagram except $0$ and $s\in \mathbb{Z}_{>0}$.

There is a conjecture that a KR module $W^{r,s}$ admits a crystal base~\cite{HKO+99,HKO+02b}, which is called a Kirillov-Reshetikhin (KR) crystal and denoted by $B^{r,s}$.
The conjecture has been solved affirmatively for all nonexceptional types~\cite{OS08}, for types $G_{2}^{(1)}$ and $D_{4}^{(3)}$~\cite{Nao18}, and for all affine types with $r$ being adjacent to $0$ or in the orbit of $0$ under a Dynkin diagram automorphism~\cite{KKM+92b}.
Recently, it has been shown that KR module $W^{r,s}$ has a crystal base whenever its classical decomposition is multiplicity free in all affine types~\cite{BS21} and has a crystal pseudobase (crystal base modulo signs) in types $E_{6,7,8}^{(1)}$, $F_{4}^{(1)}$, and $E_{6}^{(2)}$ when $r$ is near adjoint (the distance from $0$ is precisely $2$)~\cite{NS21}.

On the other hand, the notion of perfect crystal was introduced in \cite{KKM+92a}.
Hatayama et al. conjectured that KR crystals $B^{r,s}$ under some condition on $s$ are perfect crystals~\cite{HKO+02b}.
This conjecture was proved for all nonexceptional types~\cite{FOS10}, for type $G_{2}^{(1)}$ with $r=1$~\cite{Yam98} and $D_{4}^{(3)}$ with $r=1$~\cite{KMOY07}.
However, it is still an open problem for another affine types.
In this paper, we prove that KR crystals $B^{r,s}$ of types $E_{6}^{(1)}$ and $E_{7}^{(1)}$ with $r$ being minuscule nodes ($r=1,6$ for $E_{6}^{(1)}$ and $r=7$ for $E_{7}^{(1)}$), whose combinatorial structures were constructed in \cite{BS19,JS10}, are perfect crystals of level $s$.
The polytope realization of KR crystals, which was introduced in \cite{Jan21,JK19}, plays a crucial role.

This paper is organized as follows.
In Section 2 we give the necessary background on crystals.
In Section 3 we review PBW crystals and work by Salisbury-Schultze-Tingley~\cite{SST18}.
In Section 4 we summarize the facts on crystals of the quantum nilpotent subalgebras of type $E_{6,7}^{(1)}$, which are presented in \cite{Jan21}.
Section 5 is devoted to the proof of perfectness.

\section{Crystals}

\setlength{\unitlength}{10pt}

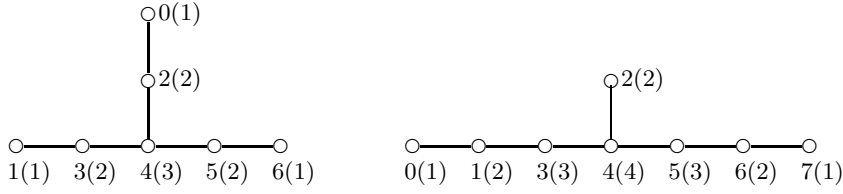
\begin{figure} 
\begin{center}
\begin{picture}(31,7)

\put(0,1.5){\circle{0.5}}
\put(2.5,1.5){\circle{0.5}}
\put(5,1.5){\circle{0.5}}
\put(7.5,1.5){\circle{0.5}}
\put(10,1.5){\circle{0.5}}
\put(5,4){\circle{0.5}}
\put(5,6.5){\circle{0.5}}

\put(0.3,1.5){\line(1,0){1.9}}
\put(2.8,1.5){\line(1,0){1.9}}
\put(5.3,1.5){\line(1,0){1.9}}
\put(7.8,1.5){\line(1,0){1.9}}
\put(5,1.8){\line(0,1){1.9}}
\put(5,4.3){\line(0,1){1.9}}

\put(0,0){\makebox(1,1){\small{$1(1)$}}}
\put(2.5,0){\makebox(1,1){\small{$3(2)$}}}
\put(5,0){\makebox(1,1){\small{$4(3)$}}}
\put(7.5,0){\makebox(1,1){\small{$5(2)$}}}
\put(10,0){\makebox(1,1){\small{$6(1)$}}}
\put(5.7,3.5){\makebox(1,1){\small{$2(2)$}}}
\put(5.7,6){\makebox(1,1){\small{$0(1)$}}}

\put(15,1.5){\circle{0.5}}
\put(17.5,1.5){\circle{0.5}}
\put(20,1.5){\circle{0.5}}
\put(22.5,1.5){\circle{0.5}}
\put(25,1.5){\circle{0.5}}
\put(27.5,1.5){\circle{0.5}}
\put(30,1.5){\circle{0.5}}
\put(22.5,4){\circle{0.5}}

\put(15.3,1.5){\line(1,0){1.9}}
\put(17.8,1.5){\line(1,0){1.9}}
\put(20.3,1.5){\line(1,0){1.9}}
\put(22.8,1.5){\line(1,0){1.9}}
\put(25.3,1.5){\line(1,0){1.9}}
\put(27.8,1.5){\line(1,0){1.9}}
\put(22.5,1.8){\line(0,1){2}}

\put(15,0){\makebox(1,1){\small{$0(1)$}}}
\put(17.5,0){\makebox(1,1){\small{$1(2)$}}}
\put(20,0){\makebox(1,1){\small{$3(3)$}}}
\put(22.5,0){\makebox(1,1){\small{$4(4)$}}}
\put(25,0){\makebox(1,1){\small{$5(3)$}}}
\put(27.5,0){\makebox(1,1){\small{$6(2)$}}}
\put(30,0){\makebox(1,1){\small{$7(1)$}}}
\put(23.2,3.5){\makebox(1,1){\small{$2(2)$}}}
\end{picture} 
\end{center}
\caption{Dynkin diagrams for $E_{6}^{(1)}$ and $E_{7}^{(1)}$ with labels from \cite{Kac90}.
Numbers in parentheses are Kac labels $a_{i}$.}
\label{fig:Dynkin}
\end{figure}

Let $\mathfrak{g}$ be an affine Lie algebra of type $E_{6,7}^{(1)}$ with Dynkin index set $I=\{0,1,\ldots,n\}$ and fundamental weights $\Lambda_{i}$ $(i\in I)$.
Let $\alpha_{i}$ and $\alpha_{i}^{\vee}$ $(i\in I)$ be the simple roots and simple coroots respectively.
Let $\delta=\sum_{i\in I}a_{i}\alpha_{i}$ be the null roots and $c^{\vee}=\sum_{i\in I}a_{i}^{\vee}\alpha_{i}^{\vee}$ the canonical central element associated to $\mathfrak{g}$, where $a_{i}$ and $a_{i}^{\vee}$ are Kac labels and dual Kac labels, which are identical in our case.
Let $P=\bigoplus_{i\in I}\mathbb{Z}\Lambda_{i}\oplus \mathbb{Z}\delta$ be the weight lattice, $P^{\vee}=\mathrm{Hom}_{\mathbb{Z}}(P,\mathbb{Z})$ the dual weight lattice, and $Q=\bigoplus_{i\in I}\mathbb{Z}\alpha_{i}$ the root lattice.
The canonical pairing $\langle \quad,\quad\rangle: P^{\vee}\times P \rightarrow \mathbb{Z}$ is defined by $\langle \alpha_{i}^{\vee},\alpha_{j}\rangle=A_{i,j}$, where $A=(A_{i,j})_{i,j\in I}$ is the generalized Cartan matrix.
We write $P_{cl}=P/\mathbb{Z}\delta$.
Let $\mathfrak{g}^{\prime}=[\mathfrak{g},\mathfrak{g}]$ be the derived subalgebra of $\mathfrak{g}$.
Let $\mathfrak{g}_{0}$ is the underlying finite-dimensional simple Lie algebra with index set $I_{0}=
I\backslash \{0\}$ and fundamental weights $\Bar{\Lambda}_{i}$ $(i\in I_{0})$.
The weight lattice is given by 
$\Bar{P}=\bigoplus_{i\in I_{0}}\mathbb{Z}\Bar{\Lambda}_{i}$ and the dual weight lattice is denoted by $\Bar{P}^{\vee}$.
Let $(\quad |\quad ): \Bar{P}\times \Bar{P}\rightarrow \mathbb{Z}$ denote the symmetric bilinear form defined by $(\alpha_{i}|\alpha_{j})=a_{i}^{\vee}a_{i}^{-1}A_{i,j}$.
Set $\theta=\sum_{i\in I_{0}}a_{i}\alpha_{i}$ the maximal root of $\mathfrak{g}_{0}$.
Note that $\alpha_{0}=-\theta \in P_{cl}$ ($a_{0}=1$).

Let $U_{q}(\mathfrak{g})$ be the quantized enveloping algebra of $\mathfrak{g}$ over $\mathbb{Q}(q)$.
A  $U_{q}(\mathfrak{g})$-crystal ($P$-weighted $I$-crystal) is a set $B$ equipped with the weight function $\mathrm{wt}: B\rightarrow P$, Kashiwara operators $e_{i},f_{i}:B\rightarrow B \cup \{\mathbf{0}\}$, and maps $\varepsilon_{i}, \varphi_{i}: B\rightarrow \mathbb{Z}\cup \{-\infty\}$ for $i\in I$ subject to the conditions:
\begin{itemize}
\item[(1)] $\varphi_{i}(b)=\varepsilon_{i}(b)+\langle \alpha_{i}^{\vee},\mathrm{wt}(b) \rangle$ for $i\in I$ and $b\in B$.
\item[(2)] if $e_{i}b\in B$, then
\begin{itemize}
\item[(a)] $\varepsilon_{i}(e_{i}b)=\varepsilon_{i}(b)-1$,
\item[(b)] $\varphi_{i}(e_{i}b)=\varphi_{i}(b)+1$, and
\item[(c)] $\mathrm{wt}(e_{i}b)=\mathrm{wt}(b)+\alpha_{i}$.
\end{itemize}
\item[(3)] if $f_{i}b\in B$, then
\begin{itemize}
\item[(a)] $\varepsilon_{i}(f_{i}b)=\varepsilon_{i}(b)+1$,
\item[(b)] $\varphi_{i}(f_{i}b)=\varphi_{i}(b)-1$, and
\item[(c)] $\mathrm{wt}(f_{i}b)=\mathrm{wt}(b)-\alpha_{i}$.
\end{itemize}
\item[(4)] $f_{i}b=b^{\prime}$ if and only if $b=e_{i}b^{\prime}$ for $b,b^{\prime}\in B$ and $i\in I$.
\item[(5)] if $\varphi_{i}(b)=-\infty$ for $b\in B$, then $e_{i}b=f_{i}b=\mathbf{0}$.
\end{itemize}

We say that the $U_{q}(\mathfrak{g})$-crystal $B$ is \emph{regular} when $\varepsilon_{i}(b)=\max \left\{ k \relmiddle| e_{i}^{k}b\neq 0 \right\}$ and $\varphi_{i}(b)=\max \left\{ k \relmiddle| f_{i}^{k}b\neq 0 \right\}$ for $b\in B$. 
Let $B$ be a regular $U_{q}(\mathfrak{g})$-crystal.
We set $\varepsilon(b)=\sum_{i\in I}\varepsilon_{i}(b)\Lambda_{i}$ and $\varphi(b)=\sum_{i\in I}\varphi_{i}(b)\Lambda_{i}$ for $b\in B$.
The weight function is given by $\mathrm{wt}(b)=\varphi(b)-\varepsilon(b)$ for $b\in B$.
We define the tensor products of $U_{q}(\mathfrak{g})$-crystals $B_{1}$ and $B_{2}$ as the crystal $B_{2}\otimes B_{1}$ that is the Cartesian product $B_{2}\times B_{1}$ with the crystal structure~\cite{BS17}
\begin{align*}
e_{i}(b_{2}\otimes b_{1})&=
\begin{cases}
e_{i}b_{2}\otimes b_{1} & \text{if } \varepsilon(b_{2})>\varphi(b_{1}), \\
b_{2}\otimes e_{i}b_{1} & \text{if } \varepsilon(b_{2})\leq\varphi(b_{1}),
\end{cases}
\\
f_{i}(b_{2}\otimes b_{1})&=
\begin{cases}
f_{i}b_{2}\otimes b_{1} & \text{if } \varepsilon(b_{2})\geq\varphi(b_{1}), \\
b_{2}\otimes f_{i}b_{1} & \text{if } \varepsilon(b_{2})<\varphi(b_{1}),
\end{cases}
\\
\varepsilon_{i}(b_{2}\otimes b_{1})&=\max (\varepsilon_{i}(b_{2}),\varepsilon_{i}(b_{1})-\langle \alpha_{i}^{\vee},\mathrm{wt}(b_{2})\rangle), \\
\varphi_{i}(b_{2}\otimes b_{1})&=\max (\varphi_{i}(b_{2}),\varphi_{i}(b_{1})+\langle \alpha_{i}^{\vee},\mathrm{wt}(b_{2})\rangle), \\
\mathrm{wt}(b_{2}\otimes b_{1})&=\mathrm{wt}(b_{2})+\mathrm{wt}(b_{1}).
\end{align*}

A $U_{q}^{\prime}(\mathfrak{g})$-crystal ($P_{cl}$-weighted $I$-crystal) and $U_{q}(\mathfrak{g}_{0})$-crystal ($\Bar{P}$-weighted $I_{0}$-crystal) are defined similarly.
For a dominant integral weight $\lambda\in \Bar{P}^{+}:=\bigoplus_{i\in I_{0}}\mathbb{Z}_{\geq 0}\Bar{\Lambda}_{i}$, we denote by $B(\lambda)$ the crystal associated to an irreducible highest weight $U_{q}(\mathfrak{g}_{0})$-module $V(\lambda)$ with highest weight $\lambda$.
For $\mu\in \Bar{P}$, let $T_{\mu}=\{t_{\mu}\}$ be a singleton crystal, where $\mathrm{wt}(t_{\mu})=\mu$ and $\varphi_{i}(t_{\mu})=-\infty$ for all $i\in I_{0}$.

\section{PBW Crystals}

The quantized enveloping algebra $U_{q}(\mathfrak{g}_{0})$ is generated by $E_{i}$, $F_{i}$, and $q^{h}$ for $i\in I_{0}$ and $h\in \Bar{P}^{\vee}$ subject to defining relations among them.
Let $U_{q}^{-}(\mathfrak{g}_{0})$ be its negative part generated by $F_{i}$'s.
For $m\in \mathbb{Z}_{>0}$, define
\[
F_{i}^{(m)}:=\frac{F_{i}^{m}}{[m]!},
\]
where $[m]!=[1][2]\cdots [m]$ and $[j]=\frac{q^{j}-q^{-j}}{q-q^{-1}}$ ($1\leq j \leq m$).
Let $A_{0}$ denote the subring of $\mathbb{Q}(q)$ of rational functions regular at $q=0$.

Let $W$ be the Weyl group of $\mathfrak{g}_{0}$ with the longest element $w_{0}$ and $\{s_{i}\}_{i\in I_{0}}$ the generating simple reflections.
Denote by $R(w_{0})$ the set of all reduced expressions for the longest element of the Weyl group.
\[
R(w_{0})=\left\{ \mathbf{i}=(i_{1},\ldots,i_{N}) \relmiddle| w_{0}=s_{i_{1}}\cdots s_{i_{N}} \right\}
\]
Elements of $R(w_{0})$ are also called as reduced words.

Given $\mathbf{i}=(i_{1},\ldots,i_{N})\in R(w_{0})$ and $\mathbf{c}=(c_{1},c_{2},\ldots,c_{N})\in \mathbb{Z}_{\geq 0}^{N}$, define
\[
b_{\mathbf{i}}(\mathbf{c})=F_{i_{1}}^{(c_{1})}T_{i_{1}}(F_{i_{2}}^{(c_{2})})\cdots T_{i_{1}}T_{i_{2}}\cdots T_{i_{N-1}}(F_{i_{N}}^{(c_{N})}),
\]
where $T_{i}$ is the Lusztig automorphism of $U_{q}(\mathfrak{g}_{0})$, which is given as $T_{i,-1}^{\prime\prime}$ in \cite{Lus93}.
Then the set $\mathbf{B}_{\mathbf{i}}:=\left\{b_{\mathbf{i}}(\mathbf{c}) \relmiddle| \mathbf{c}\in \mathbb{Z}_{\geq 0}^{N}\right\}$ is a $\mathbb{Q}(q)$-basis of $U_{q}^{-}(\mathfrak{g}_{0})$, which is called the \emph{PBW basis}.
The $A_{0}$-lattice of $U_{q}^{-}(\mathfrak{g}_{0})$ generated by $\mathbf{B}_{\mathbf{i}}$ is independent of the choice of $\mathbf{i}$, which is denoted by $L(\infty)$.
If $\pi: L(\infty)\rightarrow L(\infty)/qL(\infty)$ is the canonical projection, then the induced crystal $\pi(\mathbf{B}_{\mathbf{i}})$ is a $\mathbb{Q}$-basis of $L(\infty)/qL(\infty)$ and also independent of the choice of $\mathbf{i}$, which is denoted by  $B(\infty)$.
We identify $\mathbf{B}_{\mathbf{i}}:=\mathbb{Z}_{\geq 0}^{N}$ with a crystal $\pi(\mathbf{B}_{\mathbf{i}})$ under the map $\mathbf{c}\mapsto b_{\mathbf{i}}(\mathbf{c})$, which we call the \emph{PBW crystal}.

Let $w\in W$ be given with length $M$.
One may assume that there exists $\mathbf{i}=(i_{1},\ldots,i_{N})\in R(w_{0})$ such that $w=s_{i_{1}}\cdots s_{i_{M}}$.
The $\mathbb{Q}(q)$-subspace of $U_{q}^{-}(\mathfrak{g}_{0})$ spanned by $\mathbf{b}_{\mathbf{i}}(\mathbf{c})$ for $\mathbf{c}\in \mathbf{B}_{\mathbf{i}}$ with $c_{k}=0$ for $M+1\leq k\leq N$ is the $\mathbb{Q}(q)$-subalgebra of $U_{q}^{-}(\mathfrak{g}_{0})$.
It is independent of the choice of the reduced word of $w$.
We call this subalgebra the \emph{quantum nilpotent subalgebra} associated to $w\in W$ and denote it by $U_{q}^{-}(w)$.

Let $\Phi^{+}$ be the set of positive roots of $\mathfrak{g}_{0}$.
\begin{df}
A total order $\prec$ on $\Phi^{+}$ is called \emph{convex} if, for all triples of roots $\beta,\beta^{\prime},\beta^{\prime\prime}$ with $\beta^{\prime}=\beta+\beta^{\prime\prime}$, we have either $\beta \prec \beta^{\prime} \prec\beta^{\prime\prime}$ or $\beta^{\prime\prime} \prec \beta^{\prime} \prec \beta$.
\end{df}

\begin{thm} \cite{Pap94}
There is a bijection between $R(w_{0})$ and convex orders on $\Phi^{+}$: if $w_{0}=s_{1}\cdots s_{i_{N}}$, then the corresponding convex order $\prec$ is
\[
\beta_{1}=\alpha_{i_{1}} \prec \beta_{2}=s_{i_{1}}\alpha_{i_{2}} \prec \cdots \prec \beta_{N}=s_{i_{1}}\cdots s_{i_{N-1}}\alpha_{i_{N}}.
\]
\end{thm}

There exists a reduced word $\mathbf{i}^{\prime}$ obtained from $\mathbf{i}$ by a 3-term braid move $(i_{k},i_{k+1},i_{k+2})\rightarrow (i_{k+1},i_{k},i_{k+1})$ with $i_{k}=i_{k+2}$ if and only if $\{ \beta_{k},\beta_{k+1},\beta_{k+2} \}$ forms the positive roots of type $A_{2}$, where the corresponding convex order $\prec^{\prime}$ is given by replacing  $\beta_{k}\prec \beta_{k+1} \prec \beta_{k+2}$ with $\beta_{k+2} \prec^{\prime} \beta_{k+1} \prec^{\prime} \beta_{k}$ after the braid move.
Also there exists a reduced word $\mathbf{i}^{\prime}$ obtained from $\mathbf{i}$ by a 2-term braid move $(i_{k},i_{k+1})\rightarrow (i_{k+1},i_{k})$ if and only if $(\beta_{k}|\beta_{k+1})=0$, where the associated convex order $\prec^{\prime}$ is given by replacing $\beta_{k} \prec \beta_{k+1}$ with $\beta_{k+1} \prec^{\prime} \beta_{k}$ after the braid move.

\begin{df} \cite[Definition 4.1]{SST18}
Let $i\in I_{0}$ be given.
Then $\mathbf{i}\in R(w_{0})$ is \emph{simply braided for $i$} if one can perform a sequence of braid moves to $\mathbf{i}$ to get a word $\mathbf{i}^{\prime}$ with $i_{1}^{\prime}=i$, and each move in the sequence is either
\begin{itemize}
\item a 2-term braid move, or
\item a 3-term braid move such that $(\gamma,\gamma^{\prime},\gamma^{\prime\prime})\rightarrow (\gamma^{\prime\prime},\gamma^{\prime},\gamma)$ with $\gamma^{\prime\prime}=\alpha_{i}$.
\end{itemize}
We call $\mathbf{i}$ \emph{simply braided} if it is simply braided for all $i\in I_{0}$. 
\end{df}

Let $\sigma=(\sigma_{1},\sigma_{2},\ldots,\sigma_{s})$ be a sequence of $+$'s and $-$'s.
We cancel out an adjacent $(+,-)$-pair in $\sigma$ and repeat this process as far as possible until we get a sequence with no $-$ place to the right of $+$.
We denote the resulting sequence by $\sigma^{\mathrm{red}}$.

Given $i\in I_{0}$, suppose that $\mathbf{i}$ is simply braided for $i\in I_{0}$.
Let 
\[
\Pi_{s}=\{\gamma_{s},\gamma_{s}^{\prime},\gamma_{s}^{\prime\prime}\}
\]
be the triple of positive roots of type $A_{2}$ with $\gamma_{s}^{\prime}=\gamma_{s}+\gamma_{s}^{\prime\prime}$ and $\gamma_{s}^{\prime\prime}=\alpha_{i}$ corresponding to the $s$-th braid move for $1\leq s\leq t$.
For $\mathbf{c}\in \mathbf{B}_{\mathbf{i}}$, let
\begin{equation} \label{eq:sigma}
\sigma_{i}(\mathbf{c})=(\underbrace{-\cdots -}_{c_{\gamma_{1}^{\prime}}} \underbrace{+\cdots +}_{c_{\gamma_{1}}}\cdots
\underbrace{-\cdots -}_{c_{\gamma_{t}^{\prime}}} \underbrace{+\cdots +}_{c_{\gamma_{t}}}) .
\end{equation}

For $\beta\in \Phi^{+}$, we denote by $\mathbf{1}_{\beta}$ the element in $\mathbf{B}_{\mathbf{i}}$ where $c_{\beta}=1$ and $c_{\gamma}=0$ for $\gamma \in \Phi^{+}\backslash \{\beta\}$.

\begin{thm} \cite[Theorem 4.6]{SST18}
Let $\mathbf{i}\in R(w_{0})$ and $i\in I_{0}$.
Suppose that $\mathbf{i}$ is simply braided for $i$.
Let $\mathbf{c}\in \mathbf{B}_{\mathbf{i}}$ be given. 
\begin{itemize}
\item[(1)] If there exists $+$ in $\sigma_{i}(\mathbf{c})^{\mathrm{red}}$ and the leftmost $+$ appears in $c_{\gamma_{s}}$, then
\[
f_{i}\mathbf{c}=\mathbf{c}-\mathbf{1}_{\gamma_{s}}+\mathbf{1}_{\gamma_{s}^{\prime}}.
\]
\item[(2)] If there exists no $+$ in $\sigma_{i}(\mathbf{c})^{\mathrm{red}}$, then
\[
f_{i}\mathbf{c}=\mathbf{c}+\mathbf{1}_{\alpha_{i}}.
\]
\end{itemize}
\label{thm:SST}
\end{thm}
The rule of Theorem~\ref{thm:SST} is called the \emph{signature rule}, which was originally called the bracketing rule in \cite{SST18}.

\section{Crystals of Quantum Nilpotent Subalgebras of Type $E_{6,7}$}

We say a node $r\in I_{0}$ is \emph{special} if there exists a diagram automorphism $\phi:I\rightarrow I$ such that $\phi(0)=r$.
We say a node $r\in I_{0}$ is \emph{minuscule} if it is special and $\mathfrak{g}$ is of dual untwisted affine type.
Let $r$ be a minuscule node.
\begin{center}
\begin{tabular}{cc} \toprule
 type & $r$ \\ \midrule
 $E_{6}$ & 1,6 \\
 $E_{7}$ & 7 \\ \bottomrule
\end{tabular}
\end{center}
Set $J_{0}=I_{0}\backslash \{r\}$ and let $w^{\mathbf{J}_{0}}$ be the minimal coset representative associated to $J_{0}$.
Let us take $w_{0}=w_{\mathbf{J}_{0}}w^{\mathbf{J}_{0}}=w^{\mathbf{J}_{0}}w_{\mathbf{J}_{0}^{*}}$, where $w_{\mathbf{J}_{0}}$ is the longest element of the subalgebra $W_{J_{0}}\subset W$ generated by $s_{j}$ for $j\in J_{0}$.
The reduced word of $w^{\mathbf{J}_{0}}$ (resp. $w_{\mathbf{J}_{0}^{*}}$) is denoted by $\mathbf{i}^{\mathbf{J}_{0}}$ (resp. $\mathbf{i}_{\mathbf{J}_{0}^{*}}$) and define the reduced word $\mathbf{i}_{0}\in R(w_{0})$ by the concatenation of $\mathbf{i}^{\mathbf{J}_{0}}$ and $\mathbf{i}_{\mathbf{J}_{0}^{*}}$:
\begin{equation} \label{eq:reduced}
\mathbf{i}_{0}=\mathbf{i}^{\mathbf{J}_{0}}\mathbf{i}_{\mathbf{J}_{0}^{*}}.
\end{equation}

\begin{prp} \cite[Proposition 3.2]{Jan21}
The reduced words $\mathbf{i}^{\mathbf{J}_{0}}$ and $\mathbf{i}_{\mathbf{J}_{0}^{*}}$ are given as follows.
\begin{itemize}
\item[(1)] type $E_{6}$ with $r=1$,

$\mathbf{i}^{\mathbf{J}_{0}}=(1,3,4,5,6,2,4,5,3,4,2,1,3,4,5,6)$,

$\mathbf{i}_{\mathbf{J}_{0}^{*}}=(5,4,3,1,2,4,3,5,4,2,5,4,3,1,5,4,3,5,4,5)$.

\item[(2)] type $E_{6}$ with $r=6$,

$\mathbf{i}^{\mathbf{J}_{0}}=(6,5,4,3,1,2,4,3,5,4,2,6,5,4,3,1)$,

$\mathbf{i}_{\mathbf{J}_{0}^{*}}=(3,4,5,6,2,4,5,3,4,2,3,4,5,6,3,4,5,3,4,3)$.

\item[(3)] type $E_{7}$ with $r=7$,

$\mathbf{i}^{\mathbf{J}_{0}}=(7,6,5,4,3,1,2,4,3,5,4,2,6,5,4,3,1,7,6,5,4,3,2,4,5,6,7)$,

$\mathbf{i}_{\mathbf{J}_{0}^{*}}$ is the reduced word \eqref{eq:reduced} of type $E_{6}$ with $r=1$ or $6$. 

\end{itemize}
\end{prp}

\begin{rem}
The reduced word $\mathbf{i}_{\mathbf{J}_{0}}$ of $w_{\mathbf{J}_{0}}$ is given by $i_{\mathbf{J}_{0}}=(i_{\mathbf{J}_{0}^{*}})^{*}$, where $1^{*}=6$, $2^{*}=2$, $3^{*}=5$, $4^{*}=4$, $5^{*}=3$, and $6^{*}=1$ if $\mathfrak{g}_{0}$ is of type $E_{6}$, and $i^{*}=i$ for $i\in I_{0}$ if  $\mathfrak{g}_{0}$ is of type $E_{7}$.
\end{rem}

\begin{rem} \cite[Remark 3.4]{Jan21}
The reduced word $\mathbf{i}^{\mathbf{J}_{0}}$ is unique up to 2-term braid moves.
\end{rem}

Let $\Phi^{+}(\mathbf{J}_{0})$ be the convex-ordered set of positive roots corresponding to the reduced word $\mathbf{i}^{\mathbf{J}_{0}}$ and set 
\[
\mathbf{B}^{\mathbf{J}_{0}}=\left\{ \mathbf{c}=(c_{\beta})\in \mathbf{B}_{\mathbf{i}_{0}} \relmiddle| c_{\beta}=0 \text{ unless } \beta \in \Phi^{+}(\mathbf{J}_{0}) \right\}.
\]
Then the subcrystal $\mathbf{B}^{\mathbf{J}_{0}}$ is the quantum nilpotent subalgebra $U_{q}(w^{\mathbf{J}_{0}})$.
The convex-ordered set $\Phi^{+}(\mathbf{J}_{0})$ is give as follows.

\begin{itemize}

\item[(1)] type $E_{6}$ with $r=1$,
\begin{align*}
\Phi^{+}(\mathbf{J}_{0})=
\{&100000 \prec 101000 \prec 101100 \prec 101110 \prec 101111 \\
\prec &111100 \prec 111110 \prec 111111 \prec 111210 \prec 111211 \\
\prec &111221 \prec 112210 \prec 112211 \prec 112221 \prec 112321 \prec 122321\},
\end{align*}
where, for example, $111210$ means $\alpha_{1}+\alpha_{2}+\alpha_{3}+2\alpha_{4}+\alpha_{5}$.

\item[(2)] type $E_{6}$ with $r=6$,
\begin{align*}
\Phi^{+}(\mathbf{J}_{0})=
\{&000001 \prec 000011 \prec 000111 \prec 001111 \prec 101111 \\
\prec &010111 \prec 011111 \prec 111111 \prec 011211 \prec 111211 \\
\prec &112211 \prec 011221 \prec 111221 \prec 112221 \prec 112321 \prec 122321\}.
\end{align*}

\item[(3)] type $E_{7}$ with $r=7$,
\begin{align*}
\Phi^{+}(\mathbf{J}_{0})=
\{&0000001\prec 0000011 \prec 0000111 \prec 0001111 \prec 0011111 \prec 1011111 \\
\prec &0101111 \prec 0111111 \prec 1111111 \prec 0112111 \prec 1112111 \prec 1122111 \\
\prec &0112211 \prec 1112211 \prec 1122211 \prec 1123211 \prec 1223211 \prec 0112221 \\
\prec &1112221 \prec 1122221 \prec 1123221 \prec 1223221 \prec 1123321 \prec 1223321 \\
\prec &1224321 \prec 1234321 \prec 2234321 \}.
\end{align*}

\end{itemize}
In all cases, the first element of $\Phi^{+}(\mathbf{J}_{0})$ is $\alpha_{r}$ and the last element is the maximum root $\theta$ of $\mathfrak{g}_{0}$.

The Kashiwara operators $f_{i}$ on $\mathbf{B}^{\mathbf{J}_{0}}$ $(i\in J_{0})$ are defined by the signature rule of Theorem~\ref{thm:SST} except that (2) is replaced by
\begin{itemize}
\item[(2)'] If there exists no $+$ in $\sigma_{i}(\mathbf{c})^{\mathrm{red}}$, then $f_{i}\mathbf{c}=\mathbf{0}$.
\end{itemize}
When $i=r$, it is known that $f_{r}\mathbf{c}=\mathbf{c}+\mathbf{1}_{\alpha_{r}}$~\cite{BZ01,Lus93}.

For simplicity, we also write
\[
(\cdots \underbrace{-\cdots -}_{c_{\gamma_{s}^{\prime}}} \underbrace{+\cdots +}_{c_{\gamma_{s}}}\cdots)=
(\cdots \underbrace{-\cdots -}_{n_{\gamma_{s}^{\prime}}} \underbrace{+\cdots +}_{n_{\gamma_{s}}}\cdots)
\]
in the notation of \eqref{eq:sigma}, where $n_{\gamma_{s}^{\prime}}$ (resp. $n_{\gamma_{s}}$) is the appearance order of $\gamma_{s}^{\prime}$ (resp. $\gamma_{s}$) in the convex-ordered set of $\Phi^{+}(\mathbf{J}_{0})$.
For example, $2$ indicates $0000011$ in $\Phi^{+}(\mathbf{J}_{0})$ for type $E_{7}$ with $r=7$.

Let $i\in J_{0}$ and $\mathbf{c}\in \mathbf{B}^{\mathbf{J}_{0}}$ be given.
The sequence \eqref{eq:sigma} is given as follows.

\begin{prp} \cite[Proposition 4.3]{Jan21}

\begin{itemize}
\item[(1)] type $E_{6}$ with $r=1$,

$\sigma_{2}(\mathbf{c})=
(\underbrace{-\cdots -}_{6} \underbrace{+\cdots +}_{\textcolor{red}{3}}
\;\underbrace{-\cdots -}_{7} \underbrace{+\cdots +}_{\textcolor{red}{4}}
\;\underbrace{-\cdots -}_{8} \underbrace{+\cdots +}_{\textcolor{red}{5}}
\;\underbrace{-\cdots -}_{16} \underbrace{+\cdots +}_{\textcolor{red}{15}})$,

$\sigma_{3}(\mathbf{c})=
(\underbrace{-\cdots -}_{2} \underbrace{+\cdots +}_{\textcolor{red}{1}}
\;\underbrace{-\cdots -}_{12} \underbrace{+\cdots +}_{\textcolor{red}{9}}
\;\underbrace{-\cdots -}_{13} \underbrace{+\cdots +}_{\textcolor{red}{10}}
\;\underbrace{-\cdots -}_{14} \underbrace{+\cdots +}_{\textcolor{red}{11}})$,

$\sigma_{4}(\mathbf{c})=
(\underbrace{-\cdots -}_{3} \underbrace{+\cdots +}_{\textcolor{red}{2}}
\;\underbrace{-\cdots -}_{9} \underbrace{+\cdots +}_{\textcolor{red}{7}}
\;\underbrace{-\cdots -}_{10} \underbrace{+\cdots +}_{\textcolor{red}{8}}
\;\underbrace{-\cdots -}_{15} \underbrace{+\cdots +}_{\textcolor{red}{14}})$,

$\sigma_{5}(\mathbf{c})=
(\underbrace{-\cdots -}_{4} \underbrace{+\cdots +}_{\textcolor{red}{3}}
\;\underbrace{-\cdots -}_{7} \underbrace{+\cdots +}_{\textcolor{red}{6}}
\;\underbrace{-\cdots -}_{11} \underbrace{+\cdots +}_{\textcolor{red}{10}}
\;\underbrace{-\cdots -}_{14} \underbrace{+\cdots +}_{\textcolor{red}{13}})$,

$\sigma_{6}(\mathbf{c})=
(\underbrace{-\cdots -}_{5} \underbrace{+\cdots +}_{\textcolor{red}{4}}
\;\underbrace{-\cdots -}_{8} \underbrace{+\cdots +}_{\textcolor{red}{7}}
\;\underbrace{-\cdots -}_{10} \underbrace{+\cdots +}_{\textcolor{red}{9}}
\;\underbrace{-\cdots -}_{13} \underbrace{+\cdots +}_{\textcolor{red}{12}})$.

\item[(2)] type $E_{6}$ with $r=6$,

$\sigma_{1}(\mathbf{c})=
(\underbrace{-\cdots -}_{5} \underbrace{+\cdots +}_{\textcolor{red}{4}}
\;\underbrace{-\cdots -}_{8} \underbrace{+\cdots +}_{\textcolor{red}{7}}
\;\underbrace{-\cdots -}_{10} \underbrace{+\cdots +}_{\textcolor{red}{9}}
\;\underbrace{-\cdots -}_{13} \underbrace{+\cdots +}_{\textcolor{red}{12}})$,

$\sigma_{2}(\mathbf{c})=
(\underbrace{-\cdots -}_{6} \underbrace{+\cdots +}_{\textcolor{red}{3}}
\;\underbrace{-\cdots -}_{7} \underbrace{+\cdots +}_{\textcolor{red}{4}}
\;\underbrace{-\cdots -}_{8} \underbrace{+\cdots +}_{\textcolor{red}{5}}
\;\underbrace{-\cdots -}_{16} \underbrace{+\cdots +}_{\textcolor{red}{15}})$,

$\sigma_{3}(\mathbf{c})=
(\underbrace{-\cdots -}_{4} \underbrace{+\cdots +}_{\textcolor{red}{3}}
\;\underbrace{-\cdots -}_{7} \underbrace{+\cdots +}_{\textcolor{red}{6}}
\;\underbrace{-\cdots -}_{11} \underbrace{+\cdots +}_{\textcolor{red}{10}}
\;\underbrace{-\cdots -}_{14} \underbrace{+\cdots +}_{\textcolor{red}{13}})$,

$\sigma_{4}(\mathbf{c})=
(\underbrace{-\cdots -}_{3} \underbrace{+\cdots +}_{\textcolor{red}{2}}
\;\underbrace{-\cdots -}_{9} \underbrace{+\cdots +}_{\textcolor{red}{7}}
\;\underbrace{-\cdots -}_{10} \underbrace{+\cdots +}_{\textcolor{red}{8}}
\;\underbrace{-\cdots -}_{15} \underbrace{+\cdots +}_{\textcolor{red}{14}})$,

$\sigma_{5}(\mathbf{c})=
(\underbrace{-\cdots -}_{2} \underbrace{+\cdots +}_{\textcolor{red}{1}}
\;\underbrace{-\cdots -}_{12} \underbrace{+\cdots +}_{\textcolor{red}{9}}
\;\underbrace{-\cdots -}_{13} \underbrace{+\cdots +}_{\textcolor{red}{10}}
\;\underbrace{-\cdots -}_{14} \underbrace{+\cdots +}_{\textcolor{red}{11}})$.

\item[(3)] type $E_{7}$ with $r=7$,

\begin{align*}
\sigma_{1}(\mathbf{c})=
(\underbrace{-\cdots -}_{6} \underbrace{+\cdots +}_{\textcolor{red}{5}}
\;&\underbrace{-\cdots -}_{9} \underbrace{+\cdots +}_{\textcolor{red}{8}}
\;\underbrace{-\cdots -}_{11} \underbrace{+\cdots +}_{\textcolor{red}{10}} \\
\;&\underbrace{-\cdots -}_{14} \underbrace{+\cdots +}_{\textcolor{red}{13}}
\;\underbrace{-\cdots -}_{19} \underbrace{+\cdots +}_{\textcolor{red}{18}}
\;\underbrace{-\cdots -}_{27} \underbrace{+\cdots +}_{\textcolor{red}{26}}),
\end{align*}

\begin{align*}
\sigma_{2}(\mathbf{c})=
(\underbrace{-\cdots -}_{7} \underbrace{+\cdots +}_{\textcolor{red}{4}}
\;&\underbrace{-\cdots -}_{8} \underbrace{+\cdots +}_{\textcolor{red}{5}}
\;\underbrace{-\cdots -}_{9} \underbrace{+\cdots +}_{\textcolor{red}{6}} \\
\;&\underbrace{-\cdots -}_{17} \underbrace{+\cdots +}_{\textcolor{red}{16}}
\;\underbrace{-\cdots -}_{22} \underbrace{+\cdots +}_{\textcolor{red}{21}}
\;\underbrace{-\cdots -}_{24} \underbrace{+\cdots +}_{\textcolor{red}{23}}),
\end{align*}

\begin{align*}
\sigma_{3}(\mathbf{c})=
(\underbrace{-\cdots -}_{5} \underbrace{+\cdots +}_{\textcolor{red}{4}}
\;&\underbrace{-\cdots -}_{8} \underbrace{+\cdots +}_{\textcolor{red}{7}}
\;\underbrace{-\cdots -}_{12} \underbrace{+\cdots +}_{\textcolor{red}{11}} \\
\;&\underbrace{-\cdots -}_{15} \underbrace{+\cdots +}_{\textcolor{red}{14}}
\;\underbrace{-\cdots -}_{20} \underbrace{+\cdots +}_{\textcolor{red}{19}}
\;\underbrace{-\cdots -}_{26} \underbrace{+\cdots +}_{\textcolor{red}{25}}),
\end{align*}

\begin{align*}
\sigma_{4}(\mathbf{c})=
(\underbrace{-\cdots -}_{4} \underbrace{+\cdots +}_{\textcolor{red}{3}}
\;&\underbrace{-\cdots -}_{10} \underbrace{+\cdots +}_{\textcolor{red}{8}}
\;\underbrace{-\cdots -}_{11} \underbrace{+\cdots +}_{\textcolor{red}{9}} \\
\;&\underbrace{-\cdots -}_{16} \underbrace{+\cdots +}_{\textcolor{red}{15}}
\;\underbrace{-\cdots -}_{21} \underbrace{+\cdots +}_{\textcolor{red}{20}}
\;\underbrace{-\cdots -}_{25} \underbrace{+\cdots +}_{\textcolor{red}{24}}),
\end{align*}

\begin{align*}
\sigma_{5}(\mathbf{c})=
(\underbrace{-\cdots -}_{3} \underbrace{+\cdots +}_{\textcolor{red}{2}}
\;&\underbrace{-\cdots -}_{13} \underbrace{+\cdots +}_{\textcolor{red}{10}}
\;\underbrace{-\cdots -}_{14} \underbrace{+\cdots +}_{\textcolor{red}{11}} \\
\;&\underbrace{-\cdots -}_{15} \underbrace{+\cdots +}_{\textcolor{red}{12}}
\;\underbrace{-\cdots -}_{23} \underbrace{+\cdots +}_{\textcolor{red}{21}}
\;\underbrace{-\cdots -}_{24} \underbrace{+\cdots +}_{\textcolor{red}{22}}),
\end{align*}

\begin{align*}
\sigma_{6}(\mathbf{c})=
(\underbrace{-\cdots -}_{2} \underbrace{+\cdots +}_{\textcolor{red}{1}}
\;&\underbrace{-\cdots -}_{18} \underbrace{+\cdots +}_{\textcolor{red}{13}}
\;\underbrace{-\cdots -}_{19} \underbrace{+\cdots +}_{\textcolor{red}{14}} \\
\;&\underbrace{-\cdots -}_{20} \underbrace{+\cdots +}_{\textcolor{red}{15}}
\;\underbrace{-\cdots -}_{21} \underbrace{+\cdots +}_{\textcolor{red}{16}}
\;\underbrace{-\cdots -}_{22} \underbrace{+\cdots +}_{\textcolor{red}{17}}).
\end{align*}

\end{itemize}
\end{prp}

Let $*$ be the $\mathbb{Q}(q)$-linear anti-automorphism of $U_{q}(\mathfrak{g}_{0})$ such that $E_{i}^{*}=E_{i}$, $F_{i}^{*}=F_{i}$, and $(q^{h})^{*}=q^{-h}$ for $i\in I_{0}$ and $h\in \Bar{P}^{\vee}$.
Then $*$ induces a bijection on $B(\infty)$.
For $i\in I_{0}$, we define $e_{i}^{*}=*\circ e_{i}\circ *$ and $f_{i}^{*}=*\circ f_{i}\circ *$ on $B(\infty)$ and define $\varepsilon_{i}^{*}(b):=\max \left\{ k \relmiddle| e_{i}^{*k}b\neq 0 \right\}$ for $b\in B(\infty)$.
Using the  $\varepsilon_{i}^{*}$-statistics, we have

\begin{prp} \cite[Proposition 4.7]{Jan21}
\begin{itemize}
\item[(1)]
$\mathbf{B}^{\mathbf{J}_{0}}=\left\{ \mathbf{c}\in \mathbf{B}_{\mathbf{i}_{0}} \relmiddle| \varepsilon_{i}^{*}(\mathbf{c})=0 \text{ for } i\in J_{0}\right\}$.
\item[(2)]
Let
\[
\mathbf{B}^{\mathbf{J}_{0},s}:=\left\{ \mathbf{c}\in \mathbf{B}^{\mathbf{J_{0}}} \relmiddle| \varepsilon_{r}^{*}(\mathbf{c})\leq s \right\}
\]
for $s\geq 1$.
As $U_{q}(\mathfrak{g}_{0})$-crystals, 
\[
\mathbf{B}^{\mathbf{J}_{0},s}\otimes T_{s\Bar{\Lambda}_{r}}\simeq B(s\Bar{\Lambda}_{r}).
\]
\end{itemize}
\end{prp}
The formula for $\varepsilon_{r}^{*}(\mathbf{c})$ is given by in terms of $\mathbf{j}_{0}$-trails~\cite{BZ01}.
We just quote the result of \cite{Jan21}:
\[
\varepsilon_{r}^{*}(\mathbf{c})=\max_{l}||\mathbf{c}||_{l},
\]
where
\[
||\mathbf{c}||_{l}=\sum_{k=1}^{M}(1-m_{k}^{(l)})c_{\beta_{k}}.
\]
In the notation of \cite{Jan21}, $m_{k}^{(l)}=d_{N-k+1}(\pi_{l})$, where 
\[
(M,N)=
\begin{cases}
(16,36) & \text{for type $E_{6}$}, \\
(27,63) & \text{for type $E_{7}$},
\end{cases}
\]
and $\pi_{l}$ denote a $\mathbf{j}_{0}$-trail.

For the case of $E_{6}$, $\mathbf{m}^{(l)}=(m_{1}^{(l)},\ldots,m_{16}^{(l)})$ are 16-tuples whose entries are $0$ and $1$.
We identify $\mathbf{m}^{(l)}$ with the set of positions of $1$ in $\mathbf{m}^{(l)}$, which are common for the cases of $r=1$ and $r=6$ and listed below.

\begin{align*}
\mathbf{m}^{(1)}&=\{12,13,14,15,16\}, &\mathbf{m}^{(2)}&=\{1,13,14,15,16\}, &\mathbf{m}^{(3)}&=\{1,9,14,15,16\}, \\
\mathbf{m}^{(4)}&=\{1,9,10,15,16\}, &\mathbf{m}^{(5)}&=\{1,2,14,15,16\}, &\mathbf{m}^{(6)}&=\{1,2,10,15,16\}, \\
\mathbf{m}^{(7)}&=\{1,2,7,15,16\}, &\mathbf{m}^{(8)}&=\{1,2,7,8,16\}, &\mathbf{m}^{(9)}&=\{1,2,3,15,16\}, \\
\mathbf{m}^{(10)}&=\{1,2,3,8,16\}, & \mathbf{m}^{(11)}&=\{1,2,3,4,16\}, &\mathbf{m}^{(12)}&=\{1,2,3,4,5\}.
\end{align*}

\setlength{\unitlength}{10pt}

\begin{figure} 
\begin{center}
\begin{picture}(25,25.5)

\put( 0,24){\framebox(1,1){\small{$1$}}}
\put( 3,24){\framebox(1,1){\small{$2$}}}
\put( 6,24){\framebox(1,1){\small{$3$}}}
\put( 9,24){\framebox(1,1){\small{$4$}}}
\put(12,24){\framebox(1,1){\small{$5$}}}
\put(15,24){\framebox(1,1){\small{$6$}}}

\put( 1,24.5){\vector(1,0){2}}
\put( 4,24.5){\vector(1,0){2}}
\put( 7,24.5){\vector(1,0){2}}
\put(10,24.5){\vector(1,0){2}}
\put(13,24.5){\vector(1,0){2}}

\put( 1.5,24.5){\makebox(1,1){\tiny{$1$}}}
\put( 4.5,24.5){\makebox(1,1){\tiny{$3$}}}
\put( 7.5,24.5){\makebox(1,1){\tiny{$4$}}}
\put(10.5,24.5){\makebox(1,1){\tiny{$5$}}}
\put(13.5,24.5){\makebox(1,1){\tiny{$6$}}}

\put( 9.5,24){\vector(0,-1){2}}
\put(12.5,24){\vector(0,-1){2}}
\put(15.5,24){\vector(0,-1){2}}

\put( 8.5,22.5){\makebox(1,1){\tiny{$2$}}}
\put(11.5,22.5){\makebox(1,1){\tiny{$2$}}}
\put(14.5,22.5){\makebox(1,1){\tiny{$2$}}}

\put( 9,21){\framebox(1,1){\small{$7$}}}
\put(12,21){\framebox(1,1){\small{$8$}}}
\put(15,21){\framebox(1,1){\small{$9$}}}

\put(10,21.5){\vector(1,0){2}}
\put(13,21.5){\vector(1,0){2}}

\put(10.5,21.5){\makebox(1,1){\tiny{$5$}}}
\put(13.5,21.5){\makebox(1,1){\tiny{$6$}}}

\put(12.5,21){\vector(0,-1){2}}
\put(15.5,21){\vector(0,-1){2}}

\put(11.5,19.5){\makebox(1,1){\tiny{$4$}}}
\put(14.5,19.5){\makebox(1,1){\tiny{$4$}}}

\put(12,18){\framebox(1,1){\small{$10$}}}
\put(15,18){\framebox(1,1){\small{$11$}}}
\put(18,18){\framebox(1,1){\small{$12$}}}

\put(13,18.5){\vector(1,0){2}}
\put(16,18.5){\vector(1,0){2}}

\put(13.5,18.5){\makebox(1,1){\tiny{$6$}}}
\put(16.5,18.5){\makebox(1,1){\tiny{$5$}}}

\put(12.5,18){\vector(0,-1){2}}
\put(15.5,18){\vector(0,-1){2}}
\put(18.5,18){\vector(0,-1){2}}

\put(11.5,16.5){\makebox(1,1){\tiny{$3$}}}
\put(14.5,16.5){\makebox(1,1){\tiny{$3$}}}
\put(17.5,16.5){\makebox(1,1){\tiny{$3$}}}

\put(12,15){\framebox(1,1){\small{$13$}}}
\put(15,15){\framebox(1,1){\small{$14$}}}
\put(18,15){\framebox(1,1){\small{$15$}}}
\put(21,15){\framebox(1,1){\small{$16$}}}
\put(24,15){\framebox(1,1){\small{$17$}}}

\put(13,15.5){\vector(1,0){2}}
\put(16,15.5){\vector(1,0){2}}
\put(19,15.5){\vector(1,0){2}}
\put(22,15.5){\vector(1,0){2}}

\put(13.5,15.5){\makebox(1,1){\tiny{$6$}}}
\put(16.5,15.5){\makebox(1,1){\tiny{$5$}}}
\put(19.5,15.5){\makebox(1,1){\tiny{$4$}}}
\put(22.5,15.5){\makebox(1,1){\tiny{$2$}}}

\put(12.5,15){\vector(0,-1){2}}
\put(15.5,15){\vector(0,-1){2}}
\put(18.5,15){\vector(0,-1){2}}
\put(21.5,15){\vector(0,-1){2}}
\put(24.5,15){\vector(0,-1){2}}

\put(11.5,13.5){\makebox(1,1){\tiny{$1$}}}
\put(14.5,13.5){\makebox(1,1){\tiny{$1$}}}
\put(17.5,13.5){\makebox(1,1){\tiny{$1$}}}
\put(20.5,13.5){\makebox(1,1){\tiny{$1$}}}
\put(23.5,13.5){\makebox(1,1){\tiny{$1$}}}

\put(12,12){\framebox(1,1){\small{$18$}}}
\put(15,12){\framebox(1,1){\small{$19$}}}
\put(18,12){\framebox(1,1){\small{$20$}}}
\put(21,12){\framebox(1,1){\small{$21$}}}
\put(24,12){\framebox(1,1){\small{$22$}}}

\put(13,12.5){\vector(1,0){2}}
\put(16,12.5){\vector(1,0){2}}
\put(19,12.5){\vector(1,0){2}}
\put(22,12.5){\vector(1,0){2}}

\put(13.5,12.5){\makebox(1,1){\tiny{$6$}}}
\put(16.5,12.5){\makebox(1,1){\tiny{$5$}}}
\put(19.5,12.5){\makebox(1,1){\tiny{$4$}}}
\put(22.5,12.5){\makebox(1,1){\tiny{$2$}}}

\put(21.5,12){\vector(0,-1){2}}
\put(24.5,12){\vector(0,-1){2}}

\put(20.5,10.5){\makebox(1,1){\tiny{$3$}}}
\put(23.5,10.5){\makebox(1,1){\tiny{$3$}}}

\put(21, 9){\framebox(1,1){\small{$23$}}}
\put(24, 9){\framebox(1,1){\small{$24$}}}

\put(22, 9.5){\vector(1,0){2}}

\put(22.5, 9.5){\makebox(1,1){\tiny{$2$}}}

\put(24, 6){\framebox(1,1){\small{$25$}}}
\put(24, 3){\framebox(1,1){\small{$26$}}}
\put(24, 0){\framebox(1,1){\small{$27$}}}

\put(24.5,9){\vector(0,-1){2}}
\put(24.5,6){\vector(0,-1){2}}
\put(24.5,3){\vector(0,-1){2}}

\put(23.5,7.5){\makebox(1,1){\tiny{$4$}}}
\put(23.5,4.5){\makebox(1,1){\tiny{$5$}}}
\put(23.5,1.5){\makebox(1,1){\tiny{$6$}}}
\end{picture} 
\end{center}
\caption{Crystal graph of $\mathbf{B}^{\mathbf{J}_{0},s=1}$ for $E_{6}$ with $r=1$.}
\label{fig:graph1}
\end{figure}
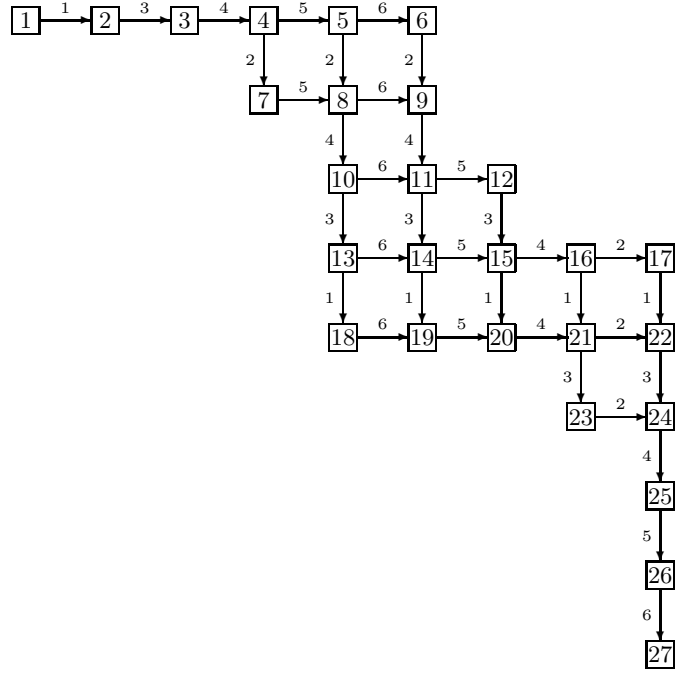 

The crystal graph of $\mathbf{B}^{\mathbf{J}_{0},s=1}$ for $E_{6}$ with $r=1$ is depicted in Figure~\ref{fig:graph1}.
The elements $\framebox{$i$}$ $(1\leq i\leq 27)$ are 16-tuples whose entries are $0$ and $1$.
We identify the element $\framebox{$i$}$ with the set of positions of $1$ in the corresponding tuple: 
$\framebox{$1$}=\{\}$,
$\framebox{$2$}=\{1\}$,
$\framebox{$3$}=\{2\}$, 
$\framebox{$4$}=\{3\}$, 
$\framebox{$5$}=\{4\}$, 
$\framebox{$6$}=\{5\}$, 
$\framebox{$7$}=\{6\}$, 
$\framebox{$8$}=\{7\}$, 
$\framebox{$9$}=\{8\}$, 
$\framebox{$10$}=\{9\}$, 
$\framebox{$11$}=\{10\}$, 
$\framebox{$12$}=\{11\}$, 
$\framebox{$13$}=\{12\}$, 
$\framebox{$14$}=\{13\}$, 
$\framebox{$15$}=\{14\}$, 
$\framebox{$16$}=\{15\}$, 
$\framebox{$17$}=\{16\}$, 
$\framebox{$18$}=\{1,12\}$, 
$\framebox{$19$}=\{1,13\}$, 
$\framebox{$20$}=\{1,14\}$, 
$\framebox{$21$}=\{1,15\}$, 
$\framebox{$22$}=\{1,16\}$, 
$\framebox{$23$}=\{2,15\}$, 
$\framebox{$24$}=\{2,16\}$, 
$\framebox{$25$}=\{3,16\}$, 
$\framebox{$26$}=\{4,16\}$, and 
$\framebox{$27$}=\{5,16\}$.
 
For the case of $E_{7}$, $\mathbf{m}^{(l)}=(m_{1}^{(l)},\ldots,m_{27}^{(l)})$ are 27-tuples whose entries are $0$ and $1$.
We identify $\mathbf{m}^{(l)}$ with the set of positions of $1$ in $\mathbf{m}^{(l)}$, which are listed below. 
 
\begin{align*}
\mathbf{m}^{(1)}&=\{18,19,20,21,23,24,25,26,27\}, &\mathbf{m}^{(2)}&=\{1,19,20,21,22,23,24,25,26,27\}, \\
\mathbf{m}^{(3)}&=\{1,13,20,21,22,23,24,25,26,27\}, &\mathbf{m}^{(4)}&=\{1,13,14,21,22,23,24,25,26,27\}, \\
\mathbf{m}^{(5)}&=\{1,13,14,15,22,23,24,25,26,27\}, &\mathbf{m}^{(6)}&=\{1,13,14,15,16,23,24,25,26,27\}, \\
\mathbf{m}^{(7)}&=\{1,2,20,21,22,23,24,25,26,27\}, &\mathbf{m}^{(8)}&=\{1,2,14,21,22,23,24,25,26,27\},
\end{align*}
\begin{align*}
\mathbf{m}^{(9)}&=\{1,2,14,15,22,23,24,25,26,27\}, &\mathbf{m}^{(10)}&=\{1,2,14,15,16,23,24,25,26,27\}, \\
\mathbf{m}^{(11)}&=\{1,2,10,21,22,23,24,25,26,27\}, &\mathbf{m}^{(12)}&=\{1,2,10,15,22,23,24,25,26,27\}, \\
\mathbf{m}^{(13)}&=\{1,2,10,15,16,23,24,25,26,27\}, &\mathbf{m}^{(14)}&=\{1,2,10,11,22,23,24,25,26,27\}, \\
\mathbf{m}^{(15)}&=\{1,2,10,11,16,23,24,25,26,27\}, &\mathbf{m}^{(16)}&=\{1,2,3,21,22,23,24,25,26,27\}, \\
\end{align*}
\begin{align*}
\mathbf{m}^{(17)}&=\{1,2,3,15,22,23,24,25,26,27\}, &\mathbf{m}^{(18)}&=\{1,2,3,15,16,23,24,25,26,27\}, \\
\mathbf{m}^{(19)}&=\{1,2,3,11,22,23,24,25,26,27\}, &\mathbf{m}^{(20)}&=\{1,2,3,8,22,23,24,25,26,27\}, \\
\mathbf{m}^{(21)}&=\{1,2,3,8,16,23,24,25,26,27\}, &\mathbf{m}^{(22)}&=\{1,2,3,11,16,23,24,25,26,27\}, \\
\mathbf{m}^{(23)}&=\{1,2,3,4,22,23,24,25,26,27\}, &\mathbf{m}^{(24)}&=\{1,2,3,4,16,23,24,25,26,27\}, \\
\end{align*}
\begin{align*}
\mathbf{m}^{(25)}&=\{1,2,3,4,9,23,24,25,26,27\}, &\mathbf{m}^{(26)}&=\{1,2,3,4,9,12,24,25,26,27\}, \\
\mathbf{m}^{(27)}&=\{1,2,3,4,9,12,21,25,26,27\}, &\mathbf{m}^{(28)}&=\{1,2,3,4,9,12,15,25,26,27\}, \\
\mathbf{m}^{(29)}&=\{1,2,3,4,9,12,15,20,26,27\}, &\mathbf{m}^{(30)}&=\{1,2,3,8,9,12,24,25,26,27\}, \\
\mathbf{m}^{(31)}&=\{1,2,3,8,9,12,21,25,26,27\}, &\mathbf{m}^{(32)}&=\{1,2,3,8,9,12,15,25,26,27\}, \\
\end{align*}
\begin{align*}
\mathbf{m}^{(33)}&=\{1,2,3,4,5,23,24,25,26,27\}, &\mathbf{m}^{(34)}&=\{1,2,3,4,5,12,24,25,26,27\}, \\
\mathbf{m}^{(35)}&=\{1,2,3,4,5,12,21,25,26,27\}, &\mathbf{m}^{(36)}&=\{1,2,3,4,5,12,15,25,26,27\}, \\
\mathbf{m}^{(37)}&=\{1,2,3,4,5,12,15,20,26,27\}, &\mathbf{m}^{(38)}&=\{1,2,3,4,5,7,24,25,26,27\}, \\
\mathbf{m}^{(39)}&=\{1,2,3,4,5,7,21,25,26,27\}, &\mathbf{m}^{(40)}&=\{1,2,3,4,5,7,15,25,26,27\}, \\
\end{align*}
\begin{align*}
\mathbf{m}^{(41)}&=\{1,2,3,4,5,7,15,20,26,27\}, &\mathbf{m}^{(42)}&=\{1,2,3,4,5,7,11,25,26,27\}, \\
\mathbf{m}^{(43)}&=\{1,2,3,4,5,7,11,20,26,27\}, &\mathbf{m}^{(44)}&=\{1,2,3,4,5,7,11,14,26,27\}, \\
\mathbf{m}^{(45)}&=\{1,2,3,4,5,7,11,14,19,27\}, &\mathbf{m}^{(46)}&=\{1,2,3,4,5,7,8,25,26,27\}, \\
\mathbf{m}^{(47)}&=\{1,2,3,4,5,7,8,20,26,27\}, &\mathbf{m}^{(48)}&=\{1,2,3,4,5,7,8,14,26,27\}, \\
\end{align*}
\begin{align*}
\mathbf{m}^{(49)}&=\{1,2,3,4,5,7,8,14,19,27\}, &\mathbf{m}^{(50)}&=\{1,2,3,4,5,7,8,10,26,27\}, \\
\mathbf{m}^{(51)}&=\{1,2,3,4,5,7,8,10,19,27\}, &\mathbf{m}^{(52)}&=\{1,2,3,4,5,7,8,10,13,27\}, \\
\mathbf{m}^{(53)}&=\{1,2,3,4,5,7,8,10,13,18\}, &\mathbf{m}^{(54)}&=\{1,2,10,11,12,22,24,25,26,27\}, \\
\mathbf{m}^{(55)}&=\{1,2,10,11,12,16,24,25,26,27\}, &\mathbf{m}^{(56)}&=\{1,2,10,11,12,16,21,25,26,27\}, \\
\end{align*}
\begin{align*}
\mathbf{m}^{(57)}&=\{1,2,3,11,12,22,24,25,26,27\}, &\mathbf{m}^{(58)}&=\{1,2,3,11,12,16,24,25,26,27\}, \\
\mathbf{m}^{(59)}&=\{1,2,3,11,12,16,21,25,26,27\}, &\mathbf{m}^{(60)}&=\{1,2,3,8,12,22,24,25,26,27\}, \\
\mathbf{m}^{(61)}&=\{1,2,3,8,12,16,24,25,26,27\}, &\mathbf{m}^{(62)}&=\{1,2,3,8,12,16,21,25,26,27\}, \\
\mathbf{m}^{(63)}&=\{1,2,3,4,12,22,24,25,26,27\}, &\mathbf{m}^{(64)}&=\{1,2,3,4,7,22,24,25,26,27\}, \\
\end{align*}
\begin{align*}
\mathbf{m}^{(65)}&=\{1,2,3,4,7,16,24,25,26,27\}, &\mathbf{m}^{(66)}&=\{1,2,3,4,7,16,21,25,26,27\}, \\
\mathbf{m}^{(67)}&=\{1,2,3,4,7,9,24,25,26,27\}, &\mathbf{m}^{(68)}&=\{1,2,3,4,7,9,21,25,26,27\}, \\
\mathbf{m}^{(69)}&=\{1,2,3,4,7,9,15,25,26,27\}, &\mathbf{m}^{(70)}&=\{1,2,3,4,7,9,15,20,26,27\}, \\
\mathbf{m}^{(71)}&=\{1,2,3,4,7,9,11,25,26,27\}, &\mathbf{m}^{(72)}&=\{1,2,3,4,7,9,11,20,26,27\}, \\
\end{align*}
\begin{align*}
\mathbf{m}^{(73)}&=\{1,2,3,4,7,9,11,14,26,27\}, &\mathbf{m}^{(74)}&=\{1,2,3,4,7,9,11,14,19,27\}, \\
\mathbf{m}^{(75)}&=\{1,2,3,4,12,16,24,25,26,27\}, &\mathbf{m}^{(76)}&=\{1,2,3,4,12,16,21,25,26,27\}. \\
\end{align*}

\setlength{\unitlength}{10pt}

\begin{figure} 
\begin{center}
\begin{picture}(31,37.5)

\put( 0,36){\framebox(1,1){\small{$1$}}}
\put( 3,36){\framebox(1,1){\small{$2$}}}
\put( 6,36){\framebox(1,1){\small{$3$}}}
\put( 9,36){\framebox(1,1){\small{$4$}}}
\put(12,36){\framebox(1,1){\small{$5$}}}
\put(15,36){\framebox(1,1){\small{$6$}}}
\put(18,36){\framebox(1,1){\small{$7$}}}

\put( 1,36.5){\vector(1,0){2}}
\put( 4,36.5){\vector(1,0){2}}
\put( 7,36.5){\vector(1,0){2}}
\put(10,36.5){\vector(1,0){2}}
\put(13,36.5){\vector(1,0){2}}
\put(16,36.5){\vector(1,0){2}}

\put( 1.5,36.5){\makebox(1,1){\tiny{$7$}}}
\put( 4.5,36.5){\makebox(1,1){\tiny{$6$}}}
\put( 7.5,36.5){\makebox(1,1){\tiny{$5$}}}
\put(10.5,36.5){\makebox(1,1){\tiny{$4$}}}
\put(13.5,36.5){\makebox(1,1){\tiny{$3$}}}
\put(16.5,36.5){\makebox(1,1){\tiny{$1$}}}

\put(12.5,36){\vector(0,-1){2}}
\put(15.5,36){\vector(0,-1){2}}
\put(18.5,36){\vector(0,-1){2}}

\put(11.5,34.5){\makebox(1,1){\tiny{$2$}}}
\put(14.5,34.5){\makebox(1,1){\tiny{$2$}}}
\put(17.5,34.5){\makebox(1,1){\tiny{$2$}}}

\put(12,33){\framebox(1,1){\small{$8$}}}
\put(15,33){\framebox(1,1){\small{$9$}}}
\put(18,33){\framebox(1,1){\small{$10$}}}

\put(13,33.5){\vector(1,0){2}}
\put(16,33.5){\vector(1,0){2}}

\put(13.5,33.5){\makebox(1,1){\tiny{$3$}}}
\put(16.5,33.5){\makebox(1,1){\tiny{$1$}}}

\put(15.5,33){\vector(0,-1){2}}
\put(18.5,33){\vector(0,-1){2}}

\put(14.5,31.5){\makebox(1,1){\tiny{$4$}}}
\put(17.5,31.5){\makebox(1,1){\tiny{$4$}}}

\put(15,30){\framebox(1,1){\small{$11$}}}
\put(18,30){\framebox(1,1){\small{$12$}}}
\put(21,30){\framebox(1,1){\small{$13$}}}

\put(16,30.5){\vector(1,0){2}}
\put(19,30.5){\vector(1,0){2}}

\put(16.5,30.5){\makebox(1,1){\tiny{$1$}}}
\put(19.5,30.5){\makebox(1,1){\tiny{$3$}}}

\put(15.5,30){\vector(0,-1){2}}
\put(18.5,30){\vector(0,-1){2}}
\put(21.5,30){\vector(0,-1){2}}

\put(14.5,28.5){\makebox(1,1){\tiny{$5$}}}
\put(17.5,28.5){\makebox(1,1){\tiny{$5$}}}
\put(20.5,28.5){\makebox(1,1){\tiny{$5$}}}

\put(15,27){\framebox(1,1){\small{$14$}}}
\put(18,27){\framebox(1,1){\small{$15$}}}
\put(21,27){\framebox(1,1){\small{$16$}}}
\put(24,27){\framebox(1,1){\small{$17$}}}
\put(27,27){\framebox(1,1){\small{$18$}}}

\put(16,27.5){\vector(1,0){2}}
\put(19,27.5){\vector(1,0){2}}
\put(22,27.5){\vector(1,0){2}}
\put(25,27.5){\vector(1,0){2}}

\put(16.5,27.5){\makebox(1,1){\tiny{$1$}}}
\put(19.5,27.5){\makebox(1,1){\tiny{$3$}}}
\put(22.5,27.5){\makebox(1,1){\tiny{$4$}}}
\put(25.5,27.5){\makebox(1,1){\tiny{$2$}}}

\put(15.5,27){\vector(0,-1){2}}
\put(18.5,27){\vector(0,-1){2}}
\put(21.5,27){\vector(0,-1){2}}
\put(24.5,27){\vector(0,-1){2}}
\put(27.5,27){\vector(0,-1){2}}

\put(14.5,25.5){\makebox(1,1){\tiny{$6$}}}
\put(17.5,25.5){\makebox(1,1){\tiny{$6$}}}
\put(20.5,25.5){\makebox(1,1){\tiny{$6$}}}
\put(23.5,25.5){\makebox(1,1){\tiny{$6$}}}
\put(26.5,25.5){\makebox(1,1){\tiny{$6$}}}

\put(15,24){\framebox(1,1){\small{$19$}}}
\put(18,24){\framebox(1,1){\small{$20$}}}
\put(21,24){\framebox(1,1){\small{$21$}}}
\put(24,24){\framebox(1,1){\small{$22$}}}
\put(27,24){\framebox(1,1){\small{$23$}}}

\put(16,24.5){\vector(1,0){2}}
\put(19,24.5){\vector(1,0){2}}
\put(22,24.5){\vector(1,0){2}}
\put(25,24.5){\vector(1,0){2}}

\put(16.5,24.5){\makebox(1,1){\tiny{$1$}}}
\put(19.5,24.5){\makebox(1,1){\tiny{$3$}}}
\put(22.5,24.5){\makebox(1,1){\tiny{$4$}}}
\put(25.5,24.5){\makebox(1,1){\tiny{$2$}}}

\put(27,21){\framebox(1,1){\small{$24$}}}
\put(30,21){\framebox(1,1){\small{$25$}}}

\put(25,24){\vector(1,-1){2}}
\put(28,24){\vector(1,-1){2}}

\put(25,23){\makebox(1,1){\tiny{$5$}}}
\put(28,22.5){\makebox(1,1){\tiny{$5$}}}

\put(28,21.5){\vector(1,0){2}}

\put(28.5,21.5){\makebox(1,1){\tiny{$2$}}}

\put(30,18){\framebox(1,1){\small{$26$}}}
\put(30,15){\framebox(1,1){\small{$27$}}}
\put(30,12){\framebox(1,1){\small{$28$}}}

\put(30.5,21){\vector(0,-1){2}}
\put(30.5,18){\vector(0,-1){2}}
\put(30.5,15){\vector(0,-1){2}}

\put(29.5,19.5){\makebox(1,1){\tiny{$4$}}}
\put(29.5,16.5){\makebox(1,1){\tiny{$3$}}}
\put(29.5,13.5){\makebox(1,1){\tiny{$1$}}}

\put(12,21){\framebox(1,1){\small{$29$}}}
\put(15,21){\framebox(1,1){\small{$30$}}}
\put(18,21){\framebox(1,1){\small{$31$}}}
\put(21,21){\framebox(1,1){\small{$32$}}}
\put(24,21){\framebox(1,1){\small{$33$}}}

\put(15,24){\vector(-1,-1){2}}
\put(18,24){\vector(-1,-1){2}}
\put(21,24){\vector(-1,-1){2}}
\put(24,24){\vector(-1,-1){2}}
\put(27,24){\vector(-1,-1){2}}

\put(13,22.5){\makebox(1,1){\tiny{$7$}}}
\put(16,22.5){\makebox(1,1){\tiny{$7$}}}
\put(19,22.5){\makebox(1,1){\tiny{$7$}}}
\put(22,22.5){\makebox(1,1){\tiny{$7$}}}
\put(26,23){\makebox(1,1){\tiny{$7$}}}

\put(13,21.5){\vector(1,0){2}}
\put(16,21.5){\vector(1,0){2}}
\put(19,21.5){\vector(1,0){2}}
\put(22,21.5){\vector(1,0){2}}

\put(13.5,21.5){\makebox(1,1){\tiny{$1$}}}
\put(16.5,21.5){\makebox(1,1){\tiny{$3$}}}
\put(19.5,21.5){\makebox(1,1){\tiny{$4$}}}
\put(22.5,21.5){\makebox(1,1){\tiny{$2$}}}

\put(24,18){\framebox(1,1){\small{$34$}}}
\put(27,18){\framebox(1,1){\small{$35$}}}

\put(22,21){\vector(1,-1){2}}
\put(25,21){\vector(1,-1){2}}

\put(22,19.5){\makebox(1,1){\tiny{$5$}}}
\put(25,20){\makebox(1,1){\tiny{$5$}}}

\put(27,21){\vector(-1,-1){2}}
\put(30,21){\vector(-1,-1){2}}
\put(30,18){\vector(-1,-1){2}}
\put(30,15){\vector(-1,-1){2}}
\put(30,12){\vector(-1,-1){2}}

\put(26,20){\makebox(1,1){\tiny{$7$}}}
\put(28.5,20){\makebox(1,1){\tiny{$7$}}}
\put(28.5,17){\makebox(1,1){\tiny{$7$}}}
\put(28.5,14){\makebox(1,1){\tiny{$7$}}}
\put(28.5,11){\makebox(1,1){\tiny{$7$}}}

\put(25,18.5){\vector(1,0){2}}

\put(25.5,18.5){\makebox(1,1){\tiny{$2$}}}

\put(27.5,18){\vector(0,-1){2}}
\put(27.5,15){\vector(0,-1){2}}
\put(27.5,12){\vector(0,-1){2}}

\put(26.5,16.5){\makebox(1,1){\tiny{$4$}}}
\put(26.5,13.5){\makebox(1,1){\tiny{$3$}}}
\put(26.5,10.5){\makebox(1,1){\tiny{$1$}}}

\put(21,15){\framebox(1,1){\small{$36$}}}
\put(24,15){\framebox(1,1){\small{$37$}}}
\put(27,15){\framebox(1,1){\small{$38$}}}

\put(24,18){\vector(-1,-1){2}}
\put(27,18){\vector(-1,-1){2}}
\put(27,15){\vector(-1,-1){2}}
\put(27,12){\vector(-1,-1){2}}
\put(27, 9){\vector(-1,-1){2}}

\put(22.5,17){\makebox(1,1){\tiny{$6$}}}
\put(25.5,17){\makebox(1,1){\tiny{$6$}}}
\put(25.5,14){\makebox(1,1){\tiny{$6$}}}
\put(25.5,11){\makebox(1,1){\tiny{$6$}}}
\put(25.5, 8){\makebox(1,1){\tiny{$6$}}}

\put(22,15.5){\vector(1,0){2}}

\put(22.5,15.5){\makebox(1,1){\tiny{$2$}}}

\put(24,12){\framebox(1,1){\small{$40$}}}
\put(27,12){\framebox(1,1){\small{$39$}}}

\put(24.5,15){\vector(0,-1){2}}
\put(24.5,12){\vector(0,-1){2}}
\put(24.5, 9){\vector(0,-1){2}}

\put(23.5,13.5){\makebox(1,1){\tiny{$4$}}}
\put(23.5,10.5){\makebox(1,1){\tiny{$3$}}}
\put(23.5, 7.5){\makebox(1,1){\tiny{$1$}}}

\put(21, 9){\framebox(1,1){\small{$43$}}}
\put(24, 9){\framebox(1,1){\small{$42$}}}
\put(27, 9){\framebox(1,1){\small{$41$}}}

\put(24,12){\vector(-1,-1){2}}
\put(24, 9){\vector(-1,-1){2}}
\put(24, 6){\vector(-1,-1){2}}

\put(22.5,11){\makebox(1,1){\tiny{$5$}}}
\put(22.5, 8){\makebox(1,1){\tiny{$5$}}}
\put(22.5, 5){\makebox(1,1){\tiny{$5$}}}

\put(15, 6){\framebox(1,1){\small{$47$}}}
\put(18, 6){\framebox(1,1){\small{$46$}}}
\put(21, 6){\framebox(1,1){\small{$45$}}}
\put(24, 6){\framebox(1,1){\small{$44$}}}

\put(21.5, 9){\vector(0,-1){2}}

\put(20.5, 7.5){\makebox(1,1){\tiny{$3$}}}

\put(21,6.5){\vector(-1,0){2}}
\put(18,6.5){\vector(-1,0){2}}

\put(19.5,6.5){\makebox(1,1){\tiny{$4$}}}
\put(16.5,6.5){\makebox(1,1){\tiny{$2$}}}

\put(15, 3){\framebox(1,1){\small{$50$}}}
\put(18, 3){\framebox(1,1){\small{$49$}}}
\put(21, 3){\framebox(1,1){\small{$48$}}}

\put(21.5, 6){\vector(0,-1){2}}
\put(18.5, 6){\vector(0,-1){2}}
\put(15.5, 6){\vector(0,-1){2}}

\put(20.5, 4.5){\makebox(1,1){\tiny{$1$}}}
\put(17.5, 4.5){\makebox(1,1){\tiny{$1$}}}
\put(14.5, 4.5){\makebox(1,1){\tiny{$1$}}}

\put(21,3.5){\vector(-1,0){2}}
\put(18,3.5){\vector(-1,0){2}}

\put(19.5,3.5){\makebox(1,1){\tiny{$4$}}}
\put(16.5,3.5){\makebox(1,1){\tiny{$2$}}}

\put( 3, 0){\framebox(1,1){\small{$56$}}}
\put( 6, 0){\framebox(1,1){\small{$55$}}}
\put( 9, 0){\framebox(1,1){\small{$54$}}}
\put(12, 0){\framebox(1,1){\small{$53$}}}
\put(15, 0){\framebox(1,1){\small{$52$}}}
\put(18, 0){\framebox(1,1){\small{$51$}}}

\put(18.5, 3){\vector(0,-1){2}}
\put(15.5, 3){\vector(0,-1){2}}

\put(17.5, 1.5){\makebox(1,1){\tiny{$3$}}}
\put(14.5, 1.5){\makebox(1,1){\tiny{$3$}}}

\put(18,0.5){\vector(-1,0){2}}
\put(15,0.5){\vector(-1,0){2}}
\put(12,0.5){\vector(-1,0){2}}
\put( 9,0.5){\vector(-1,0){2}}
\put( 6,0.5){\vector(-1,0){2}}

\put(16.5,0.5){\makebox(1,1){\tiny{$2$}}}
\put(13.5,0.5){\makebox(1,1){\tiny{$4$}}}
\put(10.5,0.5){\makebox(1,1){\tiny{$5$}}}
\put( 7.5,0.5){\makebox(1,1){\tiny{$6$}}}
\put( 4.5,0.5){\makebox(1,1){\tiny{$7$}}}

\end{picture} 
\end{center}
\caption{Crystal graph of $\mathbf{B}^{\mathbf{J}_{0},s=1}$ for $E_{7}$ with $r=7$.}
\label{fig:graph2}
\end{figure}
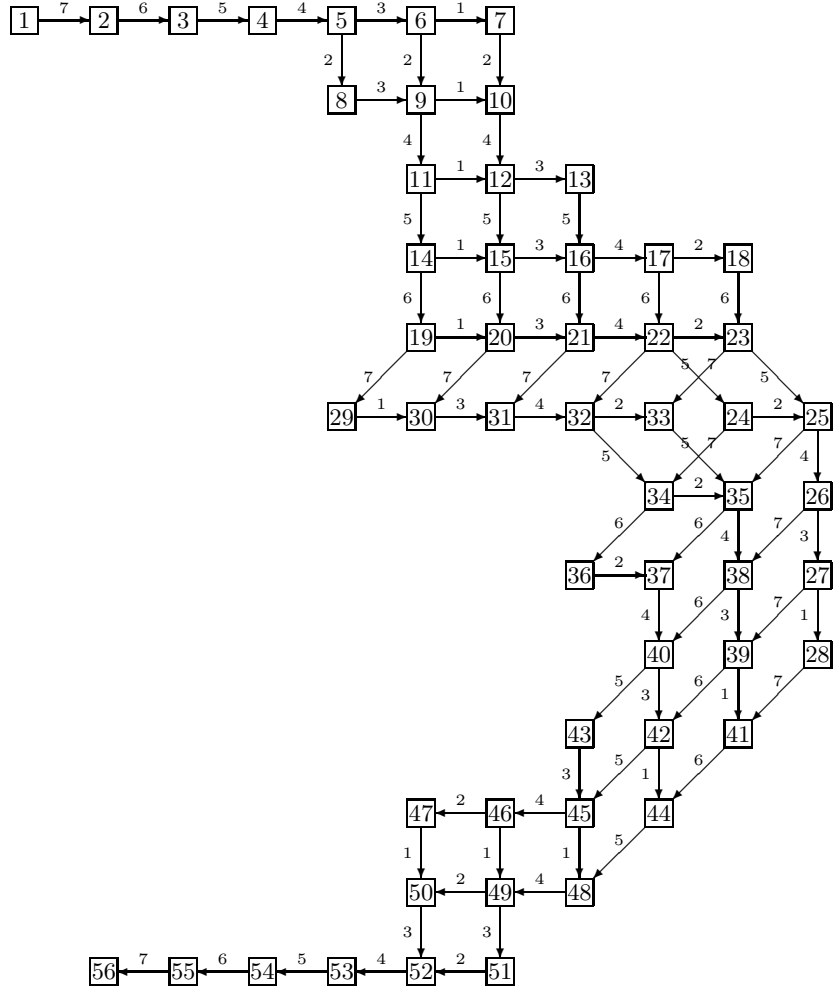

The crystal graph of $\mathbf{B}^{\mathbf{J}_{0},s=1}$ for $E_{7}$ with $r=7$ is depicted in Figure~\ref{fig:graph2}.
The elements $\framebox{$i$}$ $(1\leq i\leq 56)$ are 27-tuples whose entries are $0$ and $1$.
We identify the element $\framebox{$i$}$ with the set of positions of $1$ in the corresponding tuple:  
$\framebox{$1$}=\{\}$,
$\framebox{$2$}=\{1\}$,
$\framebox{$3$}=\{2\}$,
$\framebox{$4$}=\{3\}$,
$\framebox{$5$}=\{4\}$,
$\framebox{$6$}=\{5\}$,
$\framebox{$7$}=\{6\}$,
$\framebox{$8$}=\{7\}$,
$\framebox{$9$}=\{8\}$,
$\framebox{$10$}=\{9\}$,
$\framebox{$11$}=\{10\}$,
$\framebox{$12$}=\{11\}$,
$\framebox{$13$}=\{12\}$,
$\framebox{$14$}=\{13\}$,
$\framebox{$15$}=\{14\}$,
$\framebox{$16$}=\{15\}$,
$\framebox{$17$}=\{16\}$,
$\framebox{$18$}=\{17\}$,
$\framebox{$19$}=\{18\}$,
$\framebox{$20$}=\{19\}$,
$\framebox{$21$}=\{20\}$,
$\framebox{$22$}=\{21\}$,
$\framebox{$23$}=\{22\}$,
$\framebox{$24$}=\{23\}$,
$\framebox{$25$}=\{24\}$,
$\framebox{$26$}=\{25\}$,
$\framebox{$27$}=\{26\}$,
$\framebox{$28$}=\{27\}$,
$\framebox{$29$}=\{1,18\}$,
$\framebox{$30$}=\{1,19\}$,
$\framebox{$31$}=\{1,20\}$,
$\framebox{$32$}=\{1,21\}$,
$\framebox{$33$}=\{1,22\}$,
$\framebox{$34$}=\{1,23\}$,
$\framebox{$35$}=\{1,24\}$,
$\framebox{$36$}=\{2,23\}$,
$\framebox{$37$}=\{2,24\}$,
$\framebox{$38$}=\{1,25\}$,
$\framebox{$39$}=\{1,26\}$,
$\framebox{$40$}=\{2,25\}$,
$\framebox{$41$}=\{1,27\}$,
$\framebox{$42$}=\{2,26\}$,
$\framebox{$43$}=\{3,25\}$,
$\framebox{$44$}=\{2,27\}$,
$\framebox{$45$}=\{3,26\}$,
$\framebox{$46$}=\{4,26\}$,
$\framebox{$47$}=\{7,26\}$,
$\framebox{$48$}=\{3,27\}$,
$\framebox{$49$}=\{4,27\}$,
$\framebox{$50$}=\{7,27\}$,
$\framebox{$51$}=\{5,27\}$,
$\framebox{$52$}=\{8,27\}$,
$\framebox{$53$}=\{10,27\}$,
$\framebox{$54$}=\{13,27\}$,
$\framebox{$55$}=\{18,27\}$, and 
$\framebox{$56$}=\{1,18,27\}$.

\section{Perfectness of KR Crystals of Type $E_{6,7}^{(1)}$}

For $\mathbf{c}\in \mathbf{B}^{\mathbf{J}_{0}}$, we define
\[
f_{0}\mathbf{c}=\begin{cases}
\mathbf{c}-\mathbf{1}_{\theta} & \text{if } c_{\theta} >0, \\
\mathbf{0} & \text{otherwise,}
\end{cases}
\]
where $\theta$ is the maximum root of $\mathfrak{g}_{0}$.
Then we have
\begin{prp} \cite[Theorem 5.1]{Jan21}
For a minuscule $\Lambda_{r}$ and $s\geq 1$,
\[
\mathbf{B}^{\mathbf{J}_{0},s}\otimes T_{s\Bar{\Lambda}_{r}}\simeq B^{r,s},
\]
where $B^{r,s}$ is the KR crystal of types $E_{6}^{(1)}$ and $E_{7}^{(1)}$ associated to $s\Bar{\Lambda}_{r}$.

\end{prp}
In particular, $(0,\ldots,0)\otimes t_{s\Bar{\Lambda}_{r}} \in \mathbf{B}^{\mathbf{J}_{0},s}\otimes T_{s\Bar{\Lambda}_{r}}$ is the unique element of weight $s(\Lambda_{r}-\Lambda_{0})$.

The $0$-arrows are attached to the crystal graphs in Figures~\ref{fig:graph1} and \ref{fig:graph2}.
In the crystal graph in Figure~\ref{fig:graph1}, all the $0$-arrows are listed below.
\setlength{\unitlength}{12pt}

\begin{center}
\begin{picture}(16,4)

\put(0,0){\framebox(1,1){\small{$25$}}}
\put(3,0){\framebox(1,1){\small{$4$}}}

\put(6,0){\framebox(1,1){\small{$26$}}}
\put(9,0){\framebox(1,1){\small{$5$}}}

\put(12,0){\framebox(1,1){\small{$27$}}}
\put(15,0){\framebox(1,1){\small{$6$}}}
\put(1,0.5){\vector(1,0){2}}
\put(7,0.5){\vector(1,0){2}}
\put(13,0.5){\vector(1,0){2}}

\put(1.5,0.5){\makebox(1,1){\tiny{$0$}}}
\put(7.5,0.5){\makebox(1,1){\tiny{$0$}}}
\put(13.5,0.5){\makebox(1,1){\tiny{$0$}}}

\put(0,3){\framebox(1,1){\small{$17$}}}
\put(3,3){\framebox(1,1){\small{$1$}}}

\put(6,3){\framebox(1,1){\small{$22$}}}
\put(9,3){\framebox(1,1){\small{$2$}}}

\put(12,3){\framebox(1,1){\small{$24$}}}
\put(15,3){\framebox(1,1){\small{$3$}}}

\put(1,3.5){\vector(1,0){2}}
\put(7,3.5){\vector(1,0){2}}
\put(13,3.5){\vector(1,0){2}}

\put(1.5,3.5){\makebox(1,1){\tiny{$0$}}}
\put(7.5,3.5){\makebox(1,1){\tiny{$0$}}}
\put(13.5,3.5){\makebox(1,1){\tiny{$0$}}}

\end{picture} 
\end{center}
In the crystal graph in Figure~\ref{fig:graph2}, all the $0$-arrows are listed below.
\setlength{\unitlength}{12pt}

\begin{center}
\begin{picture}(22,7.5)

\put(0,0){\framebox(1,1){\small{$53$}}}
\put(3,0){\framebox(1,1){\small{$11$}}}

\put(6,0){\framebox(1,1){\small{$54$}}}
\put(9,0){\framebox(1,1){\small{$14$}}}

\put(12,0){\framebox(1,1){\small{$55$}}}
\put(15,0){\framebox(1,1){\small{$19$}}}

\put(18,0){\framebox(1,1){\small{$56$}}}
\put(21,0){\framebox(1,1){\small{$29$}}}

\put(1,0,5){\vector(1,0){2}}
\put(7,0,5){\vector(1,0){2}}
\put(13,0,5){\vector(1,0){2}}
\put(19,0,5){\vector(1,0){2}}

\put(1.5,0.5){\makebox(1,1){\tiny{$0$}}}
\put(7.5,0.5){\makebox(1,1){\tiny{$0$}}}
\put(13.5,0.5){\makebox(1,1){\tiny{$0$}}}
\put(19.5,0.5){\makebox(1,1){\tiny{$0$}}}

\put(0,3){\framebox(1,1){\small{$49$}}}
\put(3,3){\framebox(1,1){\small{$5$}}}

\put(6,3){\framebox(1,1){\small{$50$}}}
\put(9,3){\framebox(1,1){\small{$8$}}}

\put(12,3){\framebox(1,1){\small{$51$}}}
\put(15,3){\framebox(1,1){\small{$6$}}}

\put(18,3){\framebox(1,1){\small{$52$}}}
\put(21,3){\framebox(1,1){\small{$9$}}}

\put(1,3,5){\vector(1,0){2}}
\put(7,3,5){\vector(1,0){2}}
\put(13,3,5){\vector(1,0){2}}
\put(19,3,5){\vector(1,0){2}}

\put(1.5,3.5){\makebox(1,1){\tiny{$0$}}}
\put(7.5,3.5){\makebox(1,1){\tiny{$0$}}}
\put(13.5,3.5){\makebox(1,1){\tiny{$0$}}}
\put(19.5,3.5){\makebox(1,1){\tiny{$0$}}}

\put(0,6){\framebox(1,1){\small{$28$}}}
\put(3,6){\framebox(1,1){\small{$1$}}}

\put(6,6){\framebox(1,1){\small{$41$}}}
\put(9,6){\framebox(1,1){\small{$2$}}}

\put(12,6){\framebox(1,1){\small{$44$}}}
\put(15,6){\framebox(1,1){\small{$3$}}}

\put(18,6){\framebox(1,1){\small{$48$}}}
\put(21,6){\framebox(1,1){\small{$4$}}}

\put(1,6,5){\vector(1,0){2}}
\put(7,6,5){\vector(1,0){2}}
\put(13,6,5){\vector(1,0){2}}
\put(19,6,5){\vector(1,0){2}}

\put(1.5,6.5){\makebox(1,1){\tiny{$0$}}}
\put(7.5,6.5){\makebox(1,1){\tiny{$0$}}}
\put(13.5,6.5){\makebox(1,1){\tiny{$0$}}}
\put(19.5,6.5){\makebox(1,1){\tiny{$0$}}}
\end{picture} 
\end{center}


Let 
$P_{cl}^{+}=\left\{ \lambda\in P_{cl} \relmiddle| \langle \alpha_{i}^{\vee}, \lambda \rangle \geq 0, i\in I \right\}=\sum_{i\in I}\mathbb{Z}_{\geq 0}\Lambda_{i}$.
For $l\in \mathbb{Z}_{\geq 0}$, the set of level-$l$ weights is defined as 
$(P_{cl}^{+})_{l}=\left\{ \lambda\in P_{cl}^{+} \relmiddle| \langle c^{\vee},\lambda \rangle=l \right\}$.
Let us recall the definition of perfect crystals~\cite{KKM+92a}.
See also \cite{HK02}.

\begin{df}
For a positive integer $l>0$, a $U_{q}^{\prime}(\mathfrak{g})$-crystal $B$ is called a \emph{perfect crystal} of level $l$ if the following conditions are satisfied:
\begin{itemize}
\item[(P1)] $B$ is isomorphic to the crystal graph of a finite-dimensional $U_{q}^{\prime}(\mathfrak{g})$-module.
\item[(P2)] $B\otimes B$ is connected.
\item[(P3)] There exists a $\lambda\in P_{cl}$ such that $\mathrm{wt}(B)\subset \lambda +\sum_{i\in I_{0}}\mathbb{Z}_{\leq 0}\alpha_{i}$ and there is a unique element in $B$ of classical weight $\lambda$.
\item[(P4)] For any $b\in B$, we have $\langle c^{\vee},\varepsilon(b) \rangle\geq l$.
\item[(P5)] The maps $\varepsilon$ and $\varphi$ from $B_{\min}:=\left\{b\in B \relmiddle| \langle c^{\vee},\varepsilon(b) \rangle =l \right\}$ to $(P_{cl}^{+})_{l}$ are bijective.
\end{itemize}
\end{df}  
  
Let $B$ be KR-crystals $B^{r,s}$ $(s\geq 1)$ of type $E_{6}^{(1)}$ $(r=1,6)$ and of type $E_{7}^{(1)}$ $(r=7)$.
Since these KR crystals exist, (P1) holds.
In \cite[Corollary 6.1]{FSS07}, (P2) was proved.
The condition (P3) holds for $\lambda=s(\Lambda_{r}-\Lambda_{0})$.
In the following subsections, we show that (P4) and (P5) are satisfied.
Our main result is
\begin{thm}
Kirillov-Reshetikhin crystals $B^{r,s}$ for type $E_{6}^{(1)}$ ($r=1,6$) and $B^{7,s}$ for type $E_{7}^{(1)}$ are perfect crystals of level $s$.
\end{thm}

\subsection{$E_{6}^{(1)}$}

We give a proof for the case when $r=1$.
The proof for the case $r=6$ is almost identical.

Let $\mathbf{c}\in \mathbf{B}^{\mathbf{J}_{0},s}$.
We write $c_{\beta_{k}}=c_{k}$ and define 
$x_{l}(\mathbf{c}):=||\mathbf{c}||_{l}=\sum_{k=1}^{16}(1-m_{k}^{(l)})c_{k}$ $(l=1,2,\ldots,12)$, which must satisfy 
$x_{l}(\mathbf{c})\leq s$.
Explicitly $x_{l}(\mathbf{c})$ are written as follows.
\begin{align*}
x_{1}(\mathbf{c})&=c_{1}+c_{2}+c_{3}+c_{4}+c_{5}+c_{6}+c_{7}+c_{8}+c_{9}+c_{10}+c_{11}, \\
x_{2}(\mathbf{c})&=c_{2}+c_{3}+c_{4}+c_{5}+c_{6}+c_{7}+c_{8}+c_{9}+c_{10}+c_{11}+c_{12}, \\
x_{3}(\mathbf{c})&=c_{2}+c_{3}+c_{4}+c_{5}+c_{6}+c_{7}+c_{8}+c_{10}+c_{11}+c_{12}+c_{13}, \\
x_{4}(\mathbf{c})&=c_{2}+c_{3}+c_{4}+c_{5}+c_{6}+c_{7}+c_{8}+c_{11}+c_{12}+c_{13}+c_{14}, \\
x_{5}(\mathbf{c})&=c_{3}+c_{4}+c_{5}+c_{6}+c_{7}+c_{8}+c_{9}+c_{10}+c_{11}+c_{12}+c_{13}, \\
x_{6}(\mathbf{c})&=c_{3}+c_{4}+c_{5}+c_{6}+c_{7}+c_{8}+c_{9}+c_{11}+c_{12}+c_{13}+c_{14}, \\
x_{7}(\mathbf{c})&=c_{3}+c_{4}+c_{5}+c_{6}+c_{8}+c_{9}+c_{10}+c_{11}+c_{12}+c_{13}+c_{14}, \\
x_{8}(\mathbf{c})&=c_{3}+c_{4}+c_{5}+c_{6}+c_{9}+c_{10}+c_{11}+c_{12}+c_{13}+c_{14}+c_{15}, \\
x_{9}(\mathbf{c})&=c_{4}+c_{5}+c_{6}+c_{7}+c_{8}+c_{9}+c_{10}+c_{11}+c_{12}+c_{13}+c_{14}, \\
x_{10}(\mathbf{c})&=c_{4}+c_{5}+c_{6}+c_{7}+c_{9}+c_{10}+c_{11}+c_{12}+c_{13}+c_{14}+c_{15}, \\
x_{11}(\mathbf{c})&=c_{5}+c_{6}+c_{7}+c_{8}+c_{9}+c_{10}+c_{11}+c_{12}+c_{13}+c_{14}+c_{15}, \\
x_{12}(\mathbf{c})&=c_{6}+c_{7}+c_{8}+c_{9}+c_{10}+c_{11}+c_{12}+c_{13}+c_{14}+c_{15}+c_{16}. 
\end{align*}

Set
$c_{1}=a$, $c_{2}=b$, $c_{3}=c$, $c_{4}=d$, $c_{5}=e$, $c_{6}=f$, and
\begin{align*}
\delta_{7}&=c_{3}-c_{7}, &\delta_{8}&=c_{4}-c_{8}, &\delta_{9}&=c_{2}-c_{9}, &\delta_{10}&=c_{7}-c_{10}, \\ 
\delta_{11}&=c_{6}-c_{11}, &\delta_{12}&=c_{1}-c_{12}, &\delta_{13}&=c_{9}-c_{13}, &\delta_{14}&=c_{10}-c_{14}, \\
\delta_{15}&=c_{8}-c_{15}, &\delta_{16}&=c_{5}-c_{16}. 
\end{align*}
When $\sigma_{i}(\mathbf{c})$ for some $i$ has adjacent factors
\[
\cdots \underbrace{+\cdots +}_{k}\; \underbrace{-\cdots -}_{l}\cdots, 
\]
$\delta_{l}$ is defined as $\delta_{l}=c_{k}-c_{l}$, which is determined uniquely.

We set $s_{0}=a+2b+3c+2d+e+2f$.
Since $c_{1}$ appears only in $x_{1}(\mathbf{c})$ and $x_{1}(\mathbf{c})=s_{0}-2\delta_{7}-\delta_{8}-\delta_{9}-\delta_{10}-\delta_{11}$, 
\[
\varphi_{1}(\mathbf{c})=s-x_{1}(\mathbf{c})=s-s_{0}+2\delta_{7}+\delta_{8}+\delta_{9}+\delta_{10}+\delta_{11} 
\]
and $\varphi_{0}(\mathbf{c})=c_{16}=e-\delta_{16}$.
Set $(x)_{+}=\max (x,0)$.
From the signature rule, we have
\begin{align*}
\varphi_{2}(\mathbf{c})&=(((\delta_{7})_{+}+\delta_{8})_{+}+\delta_{16})_{+}+d-\delta_{8}-\delta_{15}, \\
\varphi_{3}(\mathbf{c})&=(((\delta_{12})_{+}+\delta_{13})_{+}+\delta_{14})_{+}+f-\delta_{11}, \\
\varphi_{4}(\mathbf{c})&=(((\delta_{9})_{+}+\delta_{10})_{+}+\delta_{15})_{+}+c-\delta_{7}-\delta_{10}-\delta_{14}, \\
\varphi_{5}(\mathbf{c})&=(((\delta_{7})_{+}+\delta_{11})_{+}+\delta_{14})_{+}+b-\delta_{9}-\delta_{13}, \text{ and} \\
\varphi_{6}(\mathbf{c})&=(((\delta_{8})_{+}+\delta_{10})_{+}+\delta_{13})_{+}+a-\delta_{12}.
\end{align*}
One may verify that $f_{i}^{\varphi_{i}(\mathbf{c})}(\mathbf{c})\in \mathbf{B}^{\mathbf{J}_{0},s}$ $(i=2,\ldots,6)$ by a case-by-case manner.
First we suppose that $\delta_{i}\leq 0$ $(7\leq i\leq 16)$.
In this case, all $\sigma_{i}(\mathbf{c})^{\mathrm{red}}$ take the form
\[
\sigma_{i}(\mathbf{c})^{\mathrm{red}}=(\cdots \underbrace{-\cdots -}_{k_{2}}\underbrace{+\cdots +}_{k_{1}}) \quad (k_{1}<k_{2})
\]
where $+$ appears only in the last factor and 
\[
\underbrace{-\cdots -}_{k_{2}}\underbrace{+\cdots +}_{k_{1}} 
\]
is the sequence of the last two factors in $\sigma_{i}(\mathbf{c})$, 
so that in $\mathbf{c}^{\prime}=f_{i}^{\varphi_{i}(\mathbf{c})}(\mathbf{c})=f_{i}^{c_{k_{1}}}(\mathbf{c})$,
$c_{k_{1}}^{\prime}=0$, $c_{k_{2}}^{\prime}=c_{k_{1}}+c_{k_{2}}$, and other components do not change.
Here, we observe that if $c_{k_{2}}$ appears in some $x_{i}(\mathbf{c})$, then $c_{k_{1}}$ always appears in the same $x_{i}(\mathbf{c})$.
Therefore, $x_{i}(\mathbf{c}^{\prime})\leq s$.
For example,
\[
\mathbf{c}^{\prime}=f_{2}^{\varphi_{2}(\mathbf{c})}(\mathbf{c})=(\cdots,0,c_{15}+c_{16})
\]
so that $x_{i}(\mathbf{c}^{\prime})=x_{i}(\mathbf{c})\leq s$ $(i=1,\ldots,7,9,12)$ and $x_{i}(\mathbf{c}^{\prime})=x_{i}(\mathbf{c})-c_{15}\leq x_{i}(\mathbf{c})\leq s$ $(i=8,10,11)$.
It remains to check that 
\begin{align*}
f_{2}^{(((\delta_{7})_{+}+\delta_{8})_{+}+\delta_{16})_{+}}(\mathbf{c})&\in \mathbf{B}^{\mathbf{J}_{0},s}, \\
f_{3}^{(((\delta_{12})_{+}+\delta_{13})_{+}+\delta_{14})_{+}}(\mathbf{c})&\in \mathbf{B}^{\mathbf{J}_{0},s}, \\
f_{4}^{(((\delta_{9})_{+}+\delta_{10})_{+}+\delta_{15})_{+}}(\mathbf{c})&\in \mathbf{B}^{\mathbf{J}_{0},s}, \\
f_{5}^{(((\delta_{7})_{+}+\delta_{11})_{+}+\delta_{14})_{+}}(\mathbf{c})&\in \mathbf{B}^{\mathbf{J}_{0},s}, \text{and} \\ 
f_{6}^{(((\delta_{8})_{+}+\delta_{10})_{+}+\delta_{13})_{+}}(\mathbf{c})&\in \mathbf{B}^{\mathbf{J}_{0},s}, 
\end{align*}
which are straightforward calculations.

Now we have
\begin{align*}
\langle c^{\vee},\varphi(\mathbf{c}) \rangle =&s -\delta_{7}-\delta_{8}-\delta_{9}-2\delta_{10}-\delta_{11}-\delta_{12}-2\delta_{13}-3\delta_{14}-2\delta_{15}-\delta_{16} \\
+&2(((\delta_{7})_{+}+\delta_{8})_{+}+\delta_{16})_{+} \\
+&2(((\delta_{12})_{+}+\delta_{13})_{+}+\delta_{14})_{+} \\
+&3(((\delta_{9})_{+}+\delta_{10})_{+}+\delta_{15})_{+} \\
+&2(((\delta_{7})_{+}+\delta_{11})_{+}+\delta_{14})_{+} \\
+&(((\delta_{8})_{+}+\delta_{10})_{+}+\delta_{13})_{+}.
\end{align*}
It is straightforward to see that  $\langle c^{\vee},\varphi(\mathbf{c}) \rangle \geq s$ so that
$\langle c^{\vee},\varphi(b) \rangle \geq s$ for $b\in B^{r,s}$. 
Since $\langle c^{\vee}, \varphi(b)-\varepsilon(b) \rangle =0$, we obtain $\langle c^{\vee}, \varepsilon(b) \rangle \geq s$, which proves (P4).

Let us write $B_{\min}=\mathbf{B}_{\min}^{\mathbf{J}_{0},s}\otimes T_{s\Bar{\Lambda}_{r}}$.
From the above arguments, it is clear that the set $\mathbf{B}_{\min}^{\mathbf{J}_{0},s}$ is a set of $\mathbf{c}\in \mathbf{B}^{\mathbf{J}_{0},s}$ such that
\begin{align*}
c_{1}&=c_{12}=a, \\
c_{2}&=c_{9}=c_{13}=b, \\
c_{3}&=c_{7}=c_{10}=c_{14}=c, \\
c_{4}&=c_{8}=c_{15}=d, \\
c_{5}&=c_{16}=e, \\
c_{6}&=c_{11}=f ,
\end{align*}
with
\[
s_{0}=a+2b+3c+2d+e+2f \leq s.
\]
For $\mathbf{c}\in \mathbf{B}_{\min}^{\mathbf{J}_{0},s}$, we have
\begin{align*}
\sigma_{2}(\mathbf{c})^{\mathrm{red}}&=
(\underbrace{-\cdots -}_{6}\;\underbrace{+\cdots +}_{15}), &
\sigma_{3}(\mathbf{c})^{\mathrm{red}}&=
(\underbrace{-\cdots -}_{2}\;\underbrace{+\cdots +}_{11}), \\
\sigma_{4}(\mathbf{c})^{\mathrm{red}}&=
(\underbrace{-\cdots -}_{3}\;\underbrace{+\cdots +}_{14}), &
\sigma_{5}(\mathbf{c})^{\mathrm{red}}&=
(\underbrace{-\cdots -}_{4}\;\underbrace{+\cdots +}_{13}), \\
\sigma_{6}(\mathbf{c})^{\mathrm{red}}&=
(\underbrace{-\cdots -}_{5}\;\underbrace{+\cdots +}_{12}),
\end{align*}
so that
\begin{align*}
\varphi_{2}(\mathbf{c})&=d, & \varphi_{3}(\mathbf{c})&=f, & \varphi_{4}(\mathbf{c})&=c, &
\varphi_{5}(\mathbf{c})=b, \\
\varphi_{6}(\mathbf{c})&=a, & \varphi_{0}(\mathbf{c})&=e, & \varphi_{1}(\mathbf{c})&=s-s_{0},
\end{align*}
and
\[
\varphi(\mathbf{c})=e\Lambda_{0}+(s-s_{0})\Lambda_{1}+d\Lambda_{2}+f\Lambda_{3}+c\Lambda_{4}+b\Lambda_{5}+a\Lambda_{6}.
\]
Let $b$ be the minimum element in $B^{r,s}$ corresponding to the minimum element $\mathbf{c} \in \mathbf{B}_{\min}^{\mathbf{J}_{0},s}$.
Then $\varphi(b)=\varphi(\mathbf{c})$ and $\mathrm{wt}(b)=\mathrm{wt}(\mathbf{c})+s(\Lambda_{r}-\Lambda_{0})$, where $\mathrm{wt}(\mathbf{c})$ is computed to be
\[
\mathrm{wt}(\mathbf{c})=e\Lambda_{0}-(s_{0}+a)\Lambda_{1}+(d-f)\Lambda_{2}-(b-f)\Lambda_{3}+(b-d)\Lambda_{5}+(a-e)\Lambda_{6}.
\]
Therefore,
\begin{align*}
\varepsilon(b)&=\varphi(b)-\mathrm{wt}(b) \\
&=(s-s_{0})\Lambda_{0}+a\Lambda_{1}+f\Lambda_{2}+b\Lambda_{3}+c\Lambda_{4}+d\Lambda_{5}+e\Lambda_{6}.
\end{align*}
Thus we have shown (P5).

\subsection{$E_{7}^{(1)}$}

Let $\mathbf{c}\in \mathbf{B}^{\mathbf{J}_{0},s}$.
Set $c_{1}=a$, $c_{2}=b$, $c_{3}=c$, $c_{4}=d$, $c_{5}=e$, $c_{6}=f$, $c_{7}=g$, and
\begin{align*}
\delta_{8}&=c_{4}-c_{8}, &\delta_{9}&=c_{5}-c_{9}, &\delta_{10}&=c_{3}-c_{10}, &\delta_{11}&=c_{8}-c_{11}, \\
\delta_{12}&=c_{7}-c_{12}, &\delta_{13}&=c_{2}-c_{13}, &\delta_{14}&=c_{10}-c_{14}, &\delta_{15}&=c_{11}-c_{15}, \\
\delta_{16}&=c_{9}-c_{16}, &\delta_{17}&=c_{6}-c_{17}, &\delta_{18}&=c_{1}-c_{18}, &\delta_{19}&=c_{13}-c_{19}, \\
\delta_{20}&=c_{14}-c_{20}, &\delta_{21}&=c_{15}-c_{21}, &\delta_{22}&=c_{16}-c_{22}, &\delta_{23}&=c_{12}-c_{23}, \\
\delta_{24}&=c_{21}-c_{24}, &\delta_{25}&=c_{20}-c_{25}, &\delta_{26}&=c_{19}-c_{26}, &\delta_{27}&=c_{18}-c_{27}.
\end{align*}
Arguments similar to the case of $E_{6}^{(1)}$ lead to $\langle c^{\vee}, \varepsilon(b) \rangle \geq s$ for $b\in B^{7,s}$, which proves (P4).

The set $\mathbf{B}_{\min}^{\mathbf{J}_{0},s}$ is a set of $\mathbf{c}\in \mathbf{B}^{\mathbf{J}_{0},s}$ such that
\begin{align*}
c_{1}&=c_{18}=c_{27}=a, \\
c_{2}&=c_{13}=c_{19}=c_{26}=b, \\
c_{3}&=c_{10}=c_{14}=c_{20}=c_{25}=c,  \\
c_{4}&=c_{8}=c_{11}=c_{15}=c_{21}=c_{24}=d, \\
c_{5}&=c_{9}=c_{16}=c_{22}=e, \\
c_{6}&=c_{17}=f, \\
c_{7}&=c_{12}=c_{23}=g, 
\end{align*}
with
\[
s_{0}=a+2b+3c+4d+3e+2f+2g \leq s.
\]
For $\mathbf{c}\in \mathbf{B}_{\min}^{\mathbf{J}_{0},s}$, we have
\begin{align*}
\sigma_{1}(\mathbf{c})^{\mathrm{red}}&=
(\underbrace{-\cdots -}_{6}\;\underbrace{+\cdots +}_{26}), &
\sigma_{2}(\mathbf{c})^{\mathrm{red}}&=
(\underbrace{-\cdots -}_{7}\;\underbrace{+\cdots +}_{23}), \\
\sigma_{3}(\mathbf{c})^{\mathrm{red}}&=
(\underbrace{-\cdots -}_{5}\;\underbrace{+\cdots +}_{25}), & 
\sigma_{4}(\mathbf{c})^{\mathrm{red}}&=
(\underbrace{-\cdots -}_{4}\;\underbrace{+\cdots +}_{24}), \\
\sigma_{5}(\mathbf{c})^{\mathrm{red}}&=
(\underbrace{-\cdots -}_{3}\;\underbrace{+\cdots +}_{22}), &
\sigma_{6}(\mathbf{c})^{\mathrm{red}}&=
(\underbrace{-\cdots -}_{2}\;\underbrace{+\cdots +}_{17}),
\end{align*}
so that
\begin{align*}
\varphi_{1}(\mathbf{c})&=b, &
\varphi_{2}(\mathbf{c})&=g, &
\varphi_{3}(\mathbf{c})&=c, &
\varphi_{4}(\mathbf{c})&=d, \\
\varphi_{5}(\mathbf{c})&=e, &
\varphi_{6}(\mathbf{c})&=f, &
\varphi_{0}(\mathbf{c})&=c_{27}=a, &
\varphi_{7}(\mathbf{c})&=s-s_{0},
\end{align*}
and
\[
\varphi(\mathbf{c})=a\Lambda_{0}+b\Lambda_{1}+g\Lambda_{2}+c\Lambda_{3}+d\Lambda_{4}+e\Lambda_{5}+f\Lambda_{6}+(s-s_{0})\Lambda_{7}.
\]
Let $b$ be the minimum element in $B^{7,s}$ corresponding to the minimum element $\mathbf{c} \in \mathbf{B}_{\min}^{\mathbf{J}_{0},s}$.
Then $\varphi(b)=\varphi(\mathbf{c})$ and
\[
\varepsilon(b)=(s-s_{0})\Lambda_{0}+f\Lambda_{1}+g\Lambda_{2}+e\Lambda_{3}+d\Lambda_{4}+c\Lambda_{5}+b\Lambda_{6}+a\Lambda_{7}.
\]
Thus we have shown (P5).

\subsection*{Acknowledgements}
The author thank Masato Okado for helpful comments.
This work was partly supported by Osaka City University Advanced Mathematical Institute (MEXT Joint Usage/Research Center on Mathematics and Theoretical Physics JPMXP0619217849).

\end{document}